\documentclass[final]{siamart}
\pdfoutput=1

\usepackage{bbm,euscript,graphicx}
\usepackage{color}
\def\cool#1{\textcolor{black}{#1}}
\def\new#1{\textcolor{black}{#1}}

\newtheorem{remark}{Remark}[section]
\def\norm#1{\|#1\|}
\def\wh#1{\widehat{#1}}
\def\wt#1{\widetilde{#1}}
\def\Re{{\rm Re}}
\def\cal{\EuScript}
\def\C{\mathbbm{C}}

\def\Cn{\C^n} 

\def\Ckk{\C^{k\times k}}

\def\Cnn{\C^{n\times n}}

\def\R{\mathbbm{R}}
\def\BA{{\bf{A}}} \def\BAs{{\bf{A}\kern-.8pt}}
\def\BB{{\bf{B}}}

\def\BD{{\bf{D}}}
\def\BE{{\bf{E}}} 
 
\def\BG{{\bf{G}}} \def\Ghat{\widehat{\BG}}
\def\BH{{\bf{H}}}
\def\BI{{\bf{I}}}
\def\BK{{\bf{K}}}

\def\BM{{\bf{M}}}
\def\BN{{\bf{N}}}
 \def\Sp{{\mbox{\textsf{p}}}}
\def\BQ{{\bf{Q}}}
\def\BR{{\bf{R}}} 
\def\Bs{{\bf{s}}} \def\BS{{\bf{S}}}  \def\Ss{{\mbox{\textsf{s}}}}
\def\BT{{\bf{T}}}
\def\BU{{\bf{U}}}  \def\Su{{\mbox{\textbf{\textsf{u}}}}}
\def\BV{{\bf{V}}} \def\BVs#1{{\bf{V}}_{\kern-1.5pt #1}} 
\def\Bv{{\bf{v}}} 
\def\Bw{{\bf{w}}} \def\Sw{{\mbox{\textbf{\textsf{w}}}}}
\def\BX{{\bf{X}}} \def\Bx{{\bf{x}}}  \def\Sx{{\mbox{\textbf{\textsf{x}}}}}
 \def\By{{\bf{y}}} 
\def\BZ{{\bf{Z}}} \def\Bz{{\bf{z}}}

\def\BLambda{\mbox{\boldmath$\Lambda$}}
\def\BDelta{\mbox{\boldmath$\Delta$}}

\def\BXi{\mbox{\boldmath$\Xi$}}
\def\Bzero{{\bf{0}}}

\def\Amu{\BA_\mu}
\def\Emu{\BE_\mu}

\def\Ker{{\rm Ker}}
\def\Ran{{\rm Ran}}
\def\SPEC{\sigma}
\def\eps{\varepsilon}
\def\PSA{\SPEC_\eps}
\def\NR{W}
\def\psabs{\alpha_\eps}
\def\nabs{\omega}
\def\dop{{\rm d}}
\def\eop{{\rm e}}
\def\iop{{\rm i}}

\mathcode`\@="8000 {\catcode`\@=\active \gdef@{\mkern1mu}}

\def\mydate{\number\day\ {\ifcase\month \or January\or February\or
              March\or April\or May\or June\or July\or August\or
              September\or October\or November\or December\fi}
\number\year}

\title{Pseudospectra of Matrix Pencils for Transient Analysis of 
         Differential-Algebraic Equations\footnotemark[1]}
\author{Mark Embree\footnotemark[2] and Blake Keeler\footnotemark[3]}

\begin{document}
\maketitle
\renewcommand{\thefootnote}{\fnsymbol{footnote}}
\footnotetext[1]{Supported 
         in part by Department of Energy grant DE-FG03-02ER25531
         and National Science Foundation grant DMS-CAREER-0449973.}
\footnotetext[2]{Department of Mathematics
and 
Computational Modeling and Data Analytics Division, 
Academy of Integrated Science,
Virginia Tech, Blacksburg, VA 24061
(\email{embree@vt.edu}).}
\footnotetext[3]{Department of Mathematics, 
University of North Carolina, Chapel Hill, NC 27599\\
(\email{bkeeler@live.unc.edu}).}


\begin{abstract}
To understand the solution of a linear, time-invariant differential-algebraic equation,
one must analyze a matrix pencil $(\BA,\BE)$ with singular $\BE$.\ \ 
Even when this pencil is stable (all its finite eigenvalues fall in the
left-half plane), the solution can exhibit transient growth before its
inevitable decay.  
When the equation results from the linearization of a nonlinear system, 
this transient growth gives a mechanism that can promote nonlinear instability.
One can enrich the conventional large-scale eigenvalue 
calculation used for linear stability analysis to identify
the potential for such transient growth.
Toward this end, we introduce a new definition of the pseudospectrum 
of a matrix pencil, use it to bound transient growth,
explain how to incorporate a physically-relevant norm,
and derive approximate pseudospectra using the
invariant subspace computed in conventional linear stability analysis.
We apply these tools to several canonical test problems in
fluid mechanics, an important source of differential-algebraic equations.
\end{abstract}



\pagestyle{myheadings} \thispagestyle{plain} 
   \markboth{M. EMBREE AND B. KEELER}{PSEUDOSPECTRA OF MATRIX PENCILS}

\section{Introduction}

Consider a linear, time invariant differential-algebraic equation (DAE)
of the general form
\begin{equation}  \label{eq:dae}
   \BE@\Bx'(t) = \BA\Bx(t), 
\end{equation}
where $\BA\in\Cnn$, $\Bx(t)\in\Cn$, and the matrix $\BE\in\Cnn$ is singular.
The singularity of $\BE$ imposes an algebraic constraint that 
any solution $\Bx(t)$ must satisfy at all $t$.  
For example, in the system
\begin{equation} \label{eq:ex3da}
\left[\!\begin{array}{ccc}
  1 & 0 & 0 \\ 0 & 1 & 0 \\ 0 & 0 & 0
\end{array}\!\right]
\left[\begin{array}{c}
x'_1(t) \\ x'_2(t) \\ x'_3(t)
\end{array}\right]
= 
\left[\!\begin{array}{rrr}
 -1 & \!\!\!-10 & 0 \\ 0 & \!-1 & 0 \\ 1 & \!1 & \!1
\end{array}\!\right]
\left[\begin{array}{c}
x_1(t) \\ x_2(t) \\ x_3(t)
\end{array}\right]
\end{equation}
the third equation gives the algebraic constraint $x_1(t)+x_2(t)+x_3(t)=0$.

Substituting the usual ansatz $\Bx(t) = \eop^{\lambda t} \Bv$ 
(for fixed $\lambda \in \C$ and $\Bv\in\Cn$) into~\cref{eq:dae}
yields the generalized eigenvalue problem
\begin{equation} \label{eq:gep}
 \BA\Bv = \lambda \BE@\Bv
\end{equation}
for the matrix pencil $(\BA,\BE)$.\ \ 
It is possible that $\BA-\lambda \BE$ is singular for all $\lambda \in \C$,
in which case the matrix pencil is \emph{singular}.
We are concerned here with the more common case of 
\emph{regular} (i.e., not singular) pencils,
where $\BA-\mu \BE$ is invertible for some $\mu\in\C$.
In this case one can find a nonzero vector $\Bv\in \Ker(\BE)$ 
(the nullspace of $\BE$) with $\BE\Bv = \Bzero$ but $\BA\Bv \ne \Bzero$;
in light of~\cref{eq:gep} we associate such $\Bv$ with the 
infinite eigenvalue $\lambda=\infty$.
This infinite eigenvalue is mapped to the 
zero eigenvalue of $\BE_\mu := (\BA-\mu\BE)^{-1}\BE$.  
In the setting of equation~\cref{eq:dae}, 
the dimension of the largest Jordan block of $\BE_\mu$ 
corresponding to a zero eigenvalue is the \emph{index}
of the differential-algebraic equation; see~\cite{CM79,KM06} for
a more detailed discussion of the index.
(The examples from fluid \new{dynamics} considered in \cref{sec:examples}
have index~2.)

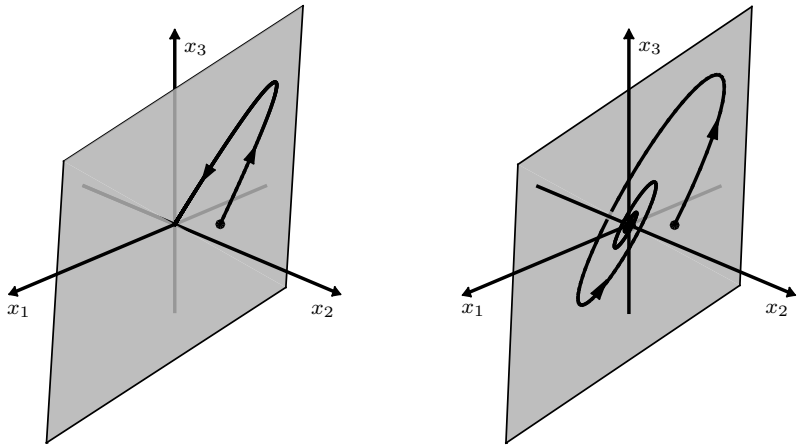
\begin{figure}[b!]
\begin{center} 
\begin{minipage}{2in}
 \includegraphics[scale=0.65]{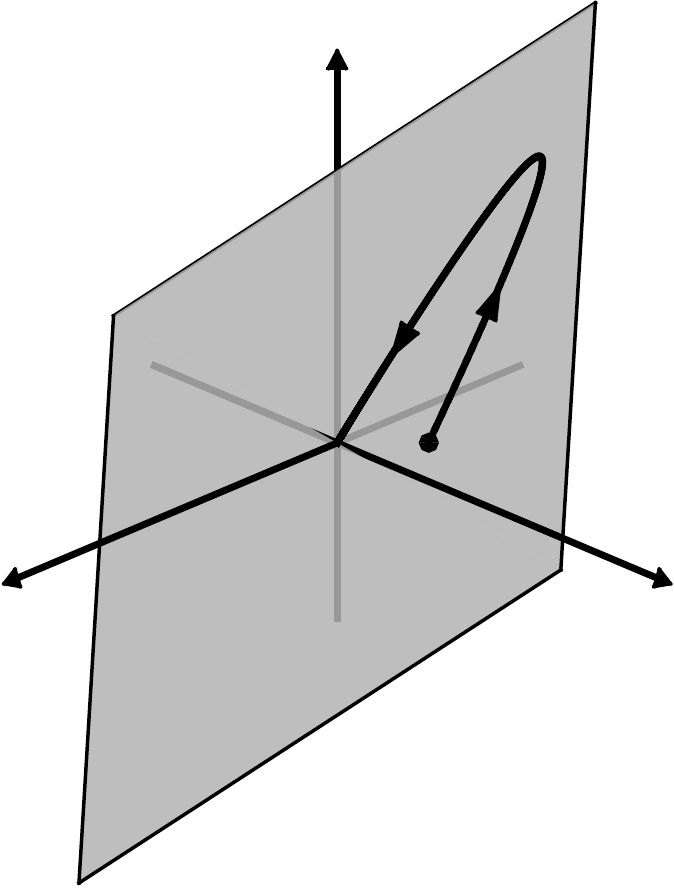}
\begin{picture}(0,0)
 \put(-130,50){\footnotesize \new{$x_1$}}
 \put(-15,50){\footnotesize \new{$x_2$}}
 \put(-63,149){\footnotesize \new{$x_3$}}
\end{picture}
\end{minipage}
\hspace*{2em}
\begin{minipage}{2in}
 \includegraphics[scale=0.65]{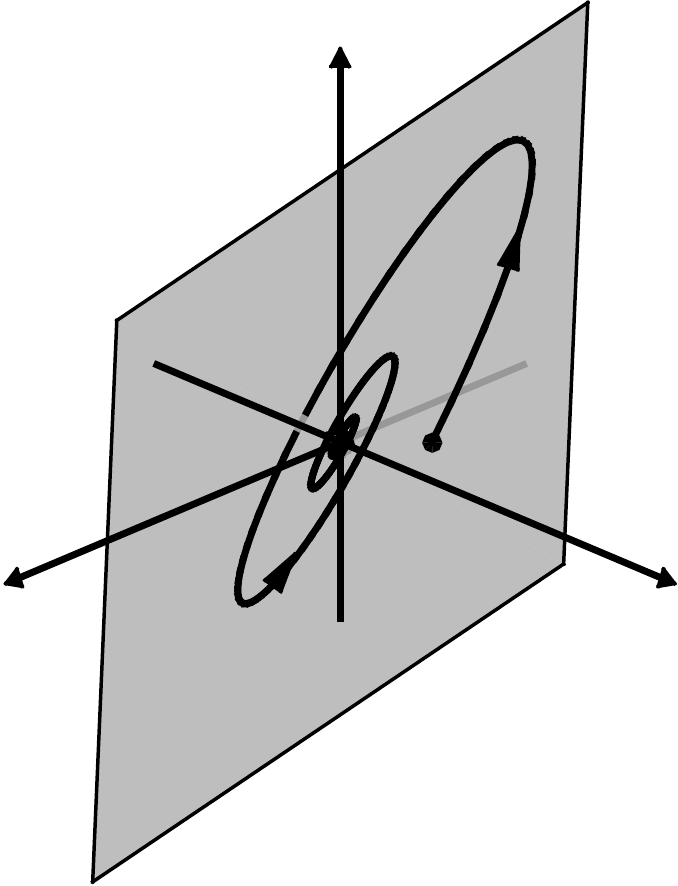}
\begin{picture}(0,0)
 \put(-130,50){\footnotesize \new{$x_1$}}
 \put(-15,50){\footnotesize \new{$x_2$}}
 \put(-63,149){\footnotesize \new{$x_3$}}
\end{picture}
\end{minipage}
\end{center}

\caption{\label{fig:ex3d}
Solutions to the DAE~$\cref{eq:ex3da}$ and the same equation with $\BA$ 
replaced by~$\cref{eq:ex3db}$, both with $\Bx(0) = [-1, 1, 0]^{\rm T}$.
The gray region indicates the plane $\{\Bx\in\R^3: x_1 +x_2+x_3=0\}$ on which
the solution is constrained to evolve.  Though both systems are asymptotically
stable, they exhibit significant transient growth:  $\|\Bx(t)\|\gg \|\Bx(0)\|$
for some values of $t>0$.
}
\end{figure}

Consider again the $(\BA,\BE)$ pair in the example~\cref{eq:ex3da}.
The pencil has spectrum $\sigma(\BA,\BE) = \{-1, \infty\}$, 
where $\lambda=-1$ has algebraic multiplicity~two.
Any initial condition must be consistent with the algebraic constraint,
i.e., $x_1(0)+x_2(0)+x_3(0)=0$, and from that initial state 
the solution will evolve in 
the two-dimensional subspace $\{\Bx\in\C^3: x_1+x_2+x_3=0\}$.
The left plot in \cref{fig:ex3d} shows the solution 
for $\Bx(0) = [-1, 1, 0]^{\rm T}$; the right plot shows
the analogous solution for the same initial condition and $\BE$,
but now with 
\begin{equation} \label{eq:ex3db}
 \BA = 
\left[\!\begin{array}{rrr}
 -1 & \!\!\!-25 & 0 \\  1 & \!-1 & 0 \\ 1 & 1 & 1
\end{array}\!\right];
\end{equation}
this modified pencil has the spectrum 
$\sigma(\BA,\BE) = \{-1+5@\iop, -1-5@\iop, \infty\}$.
In both cases the finite eigenvalues of $(\BA,\BE)$ are in the left-half plane,
so the solutions are asymptotically stable: $\Bx(t)\to \Bzero$ 
for all initial conditions that satisfy the algebraic constraint.
(For the second example the complex eigenvalues cause solutions 
to spiral toward the origin.)  
However, these examples have been designed so that
$\Bx(t)$ exhibits significant \emph{transient growth}
before eventually decaying:
there exist times $t>0$ for which $\|\Bx(t)\|\gg \|\Bx(0)\|$.
This growth is relatively modest in \cref{fig:ex3d},
compared to an increase over orders of magnitude
that could occur in some applications.

Simple eigenvalue computations alone cannot reveal the potential for transient growth,
yet such growth plays a pivotal role in dynamics.  Many DAEs of the form~\cref{eq:dae}
derive from the linear stability analysis of nonlinear dynamical systems, 
especially in fluid dynamics; see, e.g., \cite{DR81}, \cite[chap.~15]{Gun89}.  
Transient growth in the linearized system has been advanced as 
a mechanism for transition to turbulence at subcritical Reynolds numbers;
see, e.g., \cite{BDT95,BF92,Cho05,SH01,TTRD93}.
Given this possibility, classical linear stability analysis should be 
supplemented with information about transient growth, in the same way that
Gaussian elimination algorithms routinely warn when a 
matrix is severely ill-conditioned.

A variety of techniques help identify transient growth in 
the standard linear system $\Bx'(t) = \BA\Bx(t)$,
including the numerical range, pseudospectra, 
and the conditioning of a basis of eigenvectors of $\BA$;
see~\cite[Part~IV]{TE05} for a survey.
These tools do not immediately translate to the DAE setting.  
We aim to provide such a generalization,
obtaining a definition of the pseudospectrum
of a matrix pencil that preserves the algebraic structure 
of the problem, and hence is more suitable for the analysis 
of DAEs than earlier proposals in the literature.
\Cref{sec:survey} discusses these earlier definitions,
and \cref{sec:definition} describes our alternative.
This new definition is applied to derive upper and lower bounds
on the transient growth of solutions to~\cref{eq:dae} in \cref{sec:transient}.
The cost of computing pseudospectra can be a deterrent to their 
widespread adoption; thus in \cref{sec:largescale} we show
how one can readily obtain lower bounds on the proposed pseudospectra as a 
byproduct of the standard eigenvalue computation in linear stability analysis.
\Cref{sec:examples} applies these techniques to several
model problems in incompressible fluid flow, 
and \cref{sec:discrete} briefly describes how this approach 
applies to discrete-time systems with algebraic constraints.

Our primary concern here is the potential
transient growth of exact solutions of the DAE~\cref{eq:dae},
the question most relevant to linear stability analysis.  
Other definitions of pseudospectra are more appropriate when one
is concerned with \emph{uncertain} systems, as we discuss in the
next section.
We do not address other important issues, such as the challenge
of numerically generating a solution that is faithful to the
constraints~\cite{BCP96}, or understanding how the nature of the DAE 
changes under perturbations to $\BA$ and $\BE$,
which can be particularly challenging for higher index DAEs.

\section{Earlier definitions of pseudospectra of matrix pencils}
\label{sec:survey}
Throughout, we let $\sigma(\cdot)$ and $\sigma(\cdot,\cdot)$ denote the spectrum
of a matrix and matrix pencil.
For any $\eps>0$, the $\eps$-pseudospectrum $\PSA(\BA)$ 
of a matrix $\BA\in\Cnn$ is the set
\begin{align}  
 \ \ \PSA(\BA) &:= \{z\in\C: \|(z\BI-\BA)^{-1}\| > 1/\eps\} \label{eq:psa1} \\[.5em]
          &\phantom{:}= \{z\in\C: \mbox{there exists $\BDelta\in\Cnn$ with $\|\BDelta\|<\eps$ 
                and $z\in\sigma(\BA+\BDelta)$}\}, \label{eq:psa2}
\end{align}
with the convention that 
$\|\BX^{-1}\| = \infty$ when $\BX\in\Cnn$ is not invertible.
Throughout, we use the notation $\PSA(\cdot)$ with a single argument 
to denote this standard set.
The equivalence of definitions~\cref{eq:psa1} and~\cref{eq:psa2}
is fundamental to pseudospectral theory 
(see, e.g., \cite[chap.~2]{TE05} for a proof), and \new{a} cause \new{of}
ambiguity when pseudospectra are generalized
beyond the standard eigenvalue problem.
Unlike the spectrum, the pseudospectrum $\PSA(\BA)$ depends on the norm.  
\new{For now,} we let $\|\cdot\|$ denote a norm induced by an inner product,
and the associated operator norm.  (Later we will emphasize the importance 
of using physically relevant norms in our definitions.)

Since 1994 various generalizations of the $\eps$-pseudospectrum 
have been proposed for matrix pencils, e.g.,~\cite{FGNT96,HT02,LS98,Rie94,Ruh95,Dor97};
see~\cite[chap.~45]{TE05} for a comparison of these definitions.
For example, one can generalize~\cref{eq:psa1} to the pencil $(\BA,\BE)$ as
\begin{equation} \label{eq:gen1}
  \PSA(\BA,\BE) = \{z\in\C: \|(z\BE-\BA)^{-1}\| > 1/\eps\}.
\end{equation}
Alternatively, one can generalize~\cref{eq:psa2} to
\begin{align} 
\SPEC_\eps(\BA,\BE) 
          &= \{z\in\C: \mbox{there exists $\BDelta_0, \BDelta_1\in\Cnn$} 
               \label{eq:gen2} \\
          & \qquad \mbox{with $\|\BDelta_0\|<\eps@C_0$, $\|\BDelta_1\|<\eps@C_1$
          and $z\in\sigma(\BA+\BDelta_0,\BE+\BDelta_1)$}\},
               \nonumber
\end{align}
where $C_0, C_1 \ge 0$ are scaling factors that distinctly control
the size of the perturbations to each coefficient matrix.%
\footnote{The indexing of these perturbations reflects the degree 
of the coefficients $\BA$ and $\BE$ in the linear matrix pencil;
this definition further generalizes to arbitrary degree matrix 
polynomials~\cite{TH01}.}
Common choices include $C_0=C_1=1$, 
 and $C_0 = \|\BA\|$ and $C_1 = \|\BE\|$.
In fact, definition~\cref{eq:gen2} subsumes definition~\cref{eq:gen1}, since
the sets are the same when $C_0 = 1$ and $C_1 = 0$; see, e.g.,~\cite{Dor97}.

\begin{figure}[t!]
\begin{center}
\includegraphics[scale=0.3]{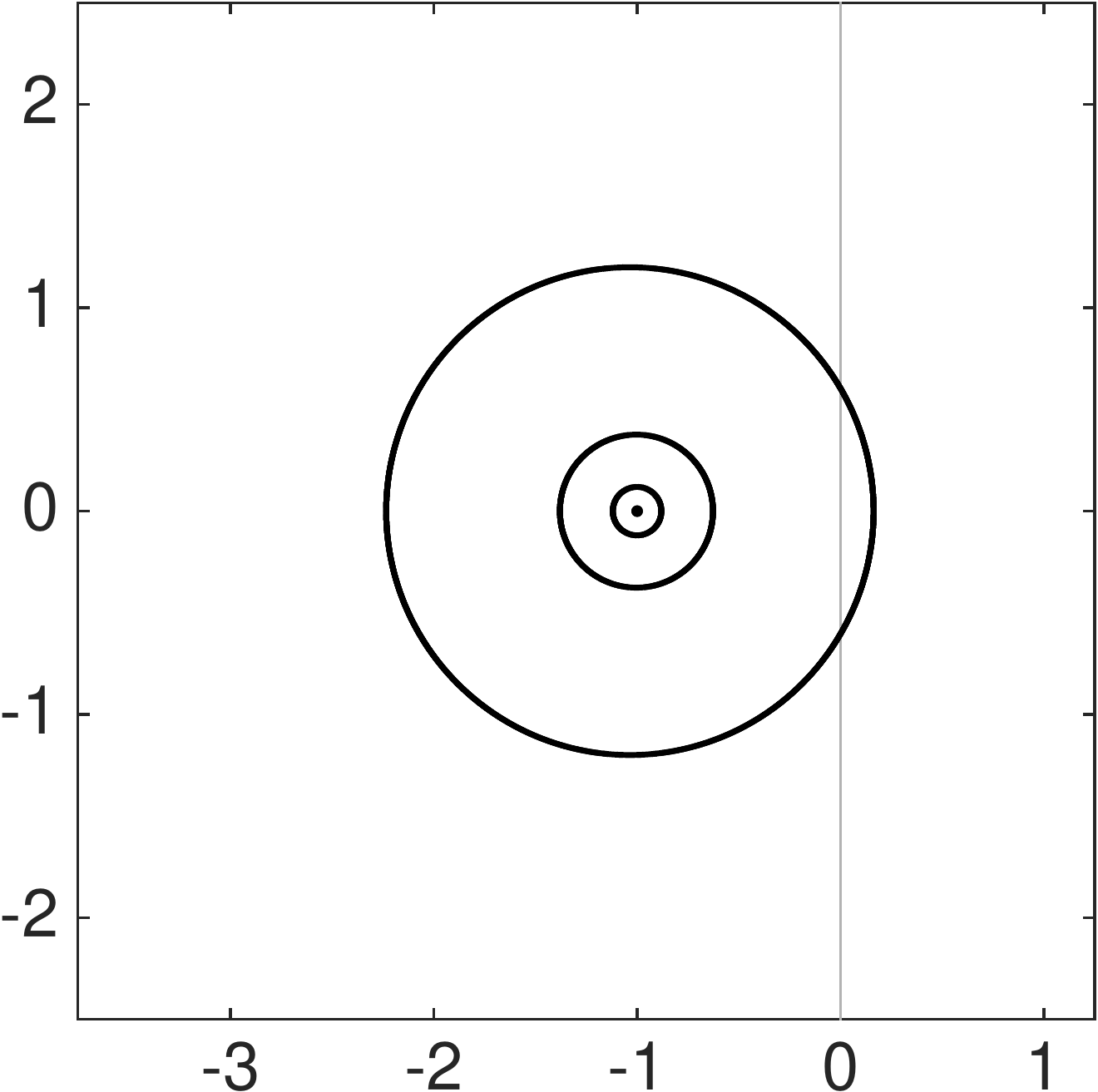}\quad
\includegraphics[scale=0.3]{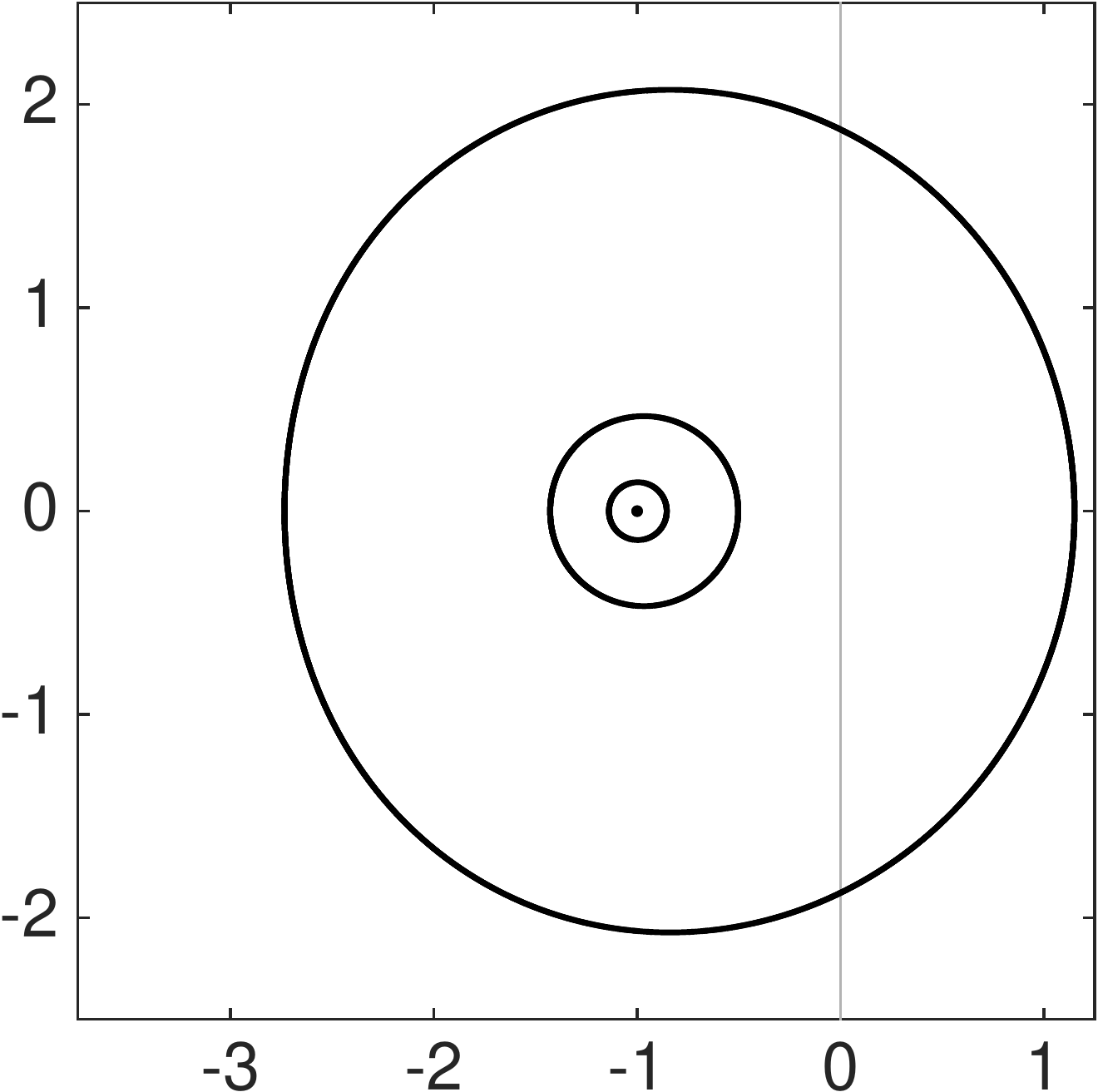}\quad
\includegraphics[scale=0.3]{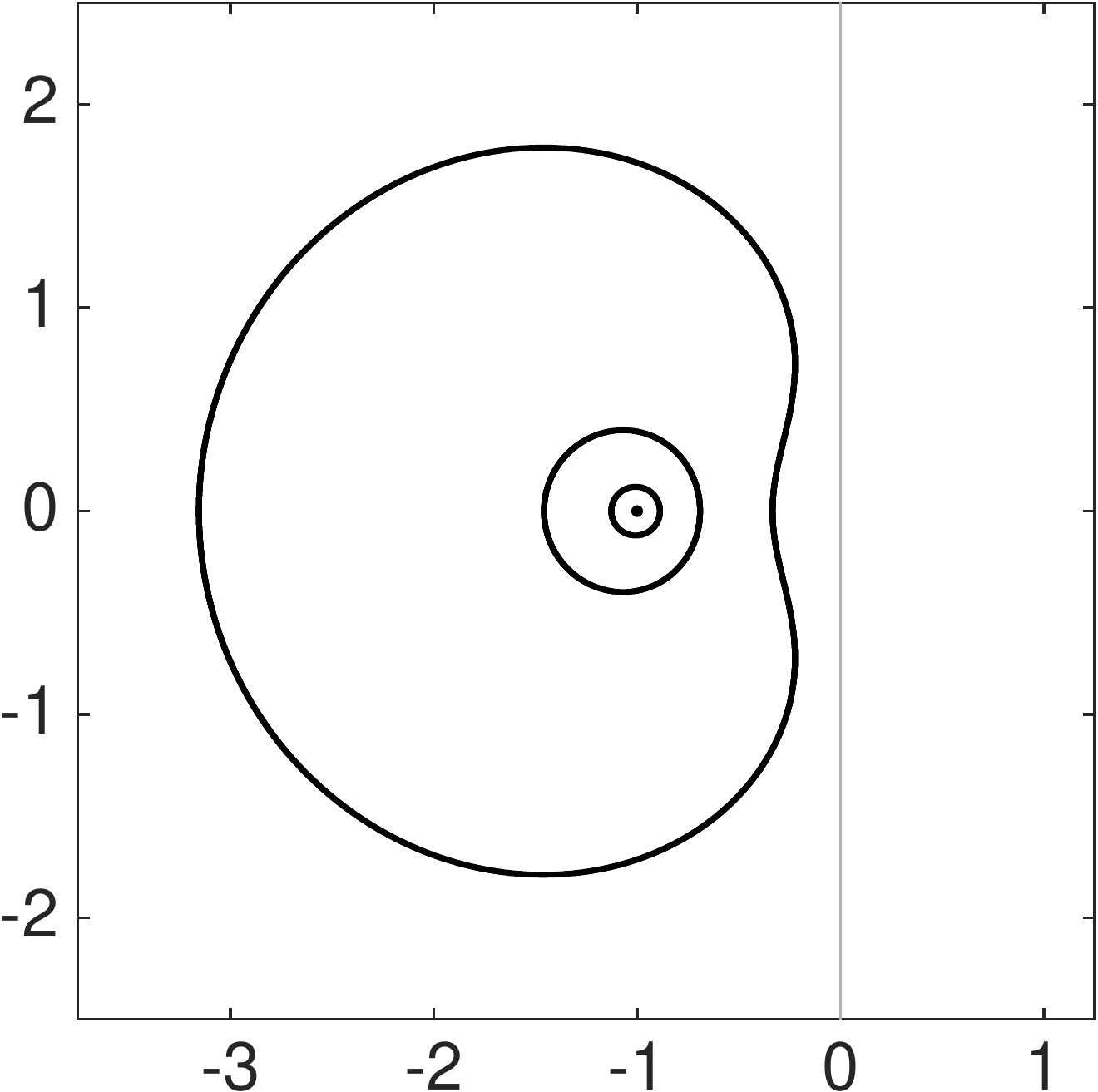}
\end{center}
\caption{\label{fig:ex1psa0}
On the left, boundaries of 
$\eps$-pseudospectra of the pencil $(\BA,\BE)$ from~$\cref{eq:ex3da}$
for $\eps= 10^{-1}$, $10^{-2}$, $10^{-3}$, 
according to definition~$\cref{eq:gen1}$, with the single
eigenvalue $\lambda=-1$.  
The middle and right plots show the same $\eps$-pseudospectra for 
$(\BT\BA,\BT\BE)$ for the two $\BT$ matrices in~$\cref{eq:T}$.
Though these pseudospectra are rather different, all three pencils
give DAEs with identical dynamics.
}
\end{figure}

Definition~\cref{eq:gen2} provides a convenient tool for assessing 
the \emph{asymptotic stability} of the solution of a DAE when the
entries of $\BA$ and $\BE$ are only known within some (bounded) uncertainty.
This definition also gives insight into the accuracy of eigenvalues 
of a matrix pencil that have been numerically computed with a 
backward-stable algorithm, and has been applied to understand the
distance of a pencil to one with a multiple eigenvalue 
(``Wilkinson's Problem'')~\cite{AAB10}.
However, as pointed out in~\cite{TE05}, this definition is unsuitable 
for analyzing the transient growth of solutions to DAEs.
To see why, premultiply equation~\cref{eq:dae} by any invertible $\BT\in\Cnn$ 
to get
\begin{equation}  \label{eq:daeS}
   \BT\BE@\Bx'(t) = \BT\BA\Bx(t).
\end{equation}
The pencil $(\BT\BA,\BT\BE)$ has the same 
spectrum as $(\BA,\BE)$ but potentially very different pseudospectra
according to definitions~\cref{eq:gen1} and~\cref{eq:gen2}.  
However, since $\BT$ has no effect on the solution $\Bx(t)$, 
it has no influence on transient dynamics.
\Cref{fig:ex1psa0} illustrates this shortcoming of
definition~\cref{eq:gen1} for the matrix pencil in equation~\cref{eq:ex3da},
comparing the $\eps$-pseudospectra of $(\BA,\BE)$ with those of 
$(\BT\BA,\BT\BE)$ for 
\begin{equation} \label{eq:T}
 \BT = \left[\begin{array}{rrr} 1 & \!\!-4 & \!16 \\ 0 & 1 & \!\!-1 \\ 0 & 0 & 1
         \end{array}\right]
   \quad \mbox{and} \quad
   \BT = \left[\begin{array}{rrr} 1 & \!\!-10 & 0 \\ 0 & 1 & 0 \\ 0 & 0 & 1
         \end{array}\right].
\end{equation}
When dealing with standard matrix pseudospectra, the rightmost extent of
$\PSA(\BA)$ in the complex plane gives crucial information about the
transient behavior of solutions to $\Bx'(t) = \BA\Bx(t)$.
(Specifics are discussed in \cref{sec:transient}.)
In each plot in \cref{fig:ex1psa0}, the outermost curve 
is the boundary of the $\eps=10^{-1}$\new{-}pseudospectrum.  
The rightmost extent of this set varies considerably 
across the three plots, even though the three
pencils define the same dynamical system, and thus 
give identical transient behavior.

To properly handle dynamics when $\BE$ is invertible, 
\cite{TE05} instead recommends Ruhe's definition~\cite{Ruh95}
\begin{align}
 \PSA(\BA,\BE) &= \PSA(\BE^{-1}\BA) \label{eq:ruhe} \\
               &= \{z\in\C: \|(z-\BE^{-1}\BA)^{-1}\| > 1/\eps\}
                \nonumber \\
               &= \{z\in\C: \|(z\BE-\BA)^{-1}\BE\| > 1/\eps\},
                \nonumber
\end{align}
emphasizing that one should use a physically
relevant norm, rather than the usual two-norm, in the definition.
(In the proper norm, \cref{eq:ruhe} can reduce to a definition
advocated by Riedel for positive definite $\BE$~\cite{Rie94}.)
Notice that the definition~\cref{eq:ruhe} is immune to the effects
of premultiplication by invertible $\BT$, since 
\[ \PSA(\BT\BA,\BT\BE) = \PSA((\BT\BE)^{-1}(\BT\BA)) = \PSA(\BE^{-1}\BA) = \PSA(\BA,\BE),\]
and, since in this case the solution of $\BE\Bx'(t)=\BA\Bx(t)$ is given by 
\begin{equation} \label{eq:etEinvA}
  \Bx(t) = \eop^{t(\BE^{-1}\BA)} \Bx(0),
\end{equation}
one can understand the transient dynamics of~\cref{eq:etEinvA} 
from standard results about the pseudospectra of $\BE^{-1}\BA$.
However, this definition is clearly insufficient for differential-algebraic
equations, where $\BE$ is not invertible.%
\footnote{For singular $\BE$, 
\cite[pp.~428--429]{TE05} tentatively
suggests a regularization approach that turns out to be 
insufficient for describing DAE dynamics.}

\section{Pseudospectra for matrix pencils derived from DAEs}
\label{sec:definition}
To begin this section let $\|\cdot\|$ denote the vector 2-norm
and the matrix norm it induces; more general norms will be 
addressed in \cref{sec:norms}.

Our definition of pseudospectra for matrix pencils derived
from DAEs follows from a simple strategy: to gain insight 
into the transient dynamics, we should base our definition on
the roles that $\BA$ and $\BE$ play in the solution formula 
for the DAE.%
\footnote{This approach amounts to defining 
pseudospectra in terms of the infinitesimal generator in the
the semigroup formula for the solution $\Bx(t)$.  
Green and Wagenknecht briefly mention the analogous definition
for delay differential equations in~\cite[sect.~4]{GW06}.
The decomposition we use here is commonly applied in
reduced order modeling for descriptor systems; see, e.g., \cite{HSS08,Sty04}.}
These solutions are typically expressed using the Drazin inverse 
(see, e.g., \cite[chap.~9]{CM79}, \cite{KM06}).
While this approach gives an algebraically elegant, compact formula,
its use of the Jordan form is computationally unappealing.
We shall essentially recapitulate the derivation 
from~\cite{CM79}, but instead use the Schur factorization.

Suppose $(\BA,\BE)$ is a regular pencil, so there exists
some $\mu\in\C$ such that $\BA-\mu\BE$ is invertible. 
For such a $\mu$ define
\[ \Amu := (\BA-\mu\BE)^{-1}\BA, \qquad
   \Emu := (\BA-\mu\BE)^{-1}\BE,
\]
and premultiply~the DAE~\cref{eq:dae} by $(\BA-\mu@\BE)^{-1}$ 
to obtain
\begin{equation} \label{eq:EmuAmu}
   \Emu \Bx'(t) = \Amu \Bx(t). 
\end{equation}
Now since 
$\Amu = (\BA-\mu\BE)^{-1}(\BA-\mu\BE+\mu\BE) = \BI + \mu@\Emu$, 
\cref{eq:EmuAmu} can be written as
\begin{equation} \label{eq:EmuEmu}
   \Emu \Bx'(t) = (\BI+\mu\Emu) \Bx(t). 
\end{equation}
Compute the Schur factorization
\begin{equation} \label{eq:schur}
 \Emu = \left[\begin{array}{cc} \BQ_\mu & \wt{\BQ}_\mu \end{array}\right]
        \left[\begin{array}{cc}  \BG_\mu & \BD_{\kern-.5pt \mu} \\ \Bzero & \BN_\mu \end{array}\right]
        \left[\begin{array}{c}  \BQ_\mu^* \\ \wt{\BQ}_\mu^*\end{array}\right],
\end{equation}
where $[\BQ_\mu\ \wt{\BQ}_\mu]\in\Cnn$ is unitary and
the diagonal of the Schur factor has been ordered so that
$\BN_\mu\in\C^{d\times d}$ is nilpotent, containing all the 
zero eigenvalues of $\Emu$, and hence $0\not\in\SPEC(\BG_\mu)$.
(This factorization can be computed via the generalized null space 
decomposition algorithm~\cite{GOS15}.
In many cases, the dimension $d$ is known directly from the
application, as is common in fluid mechanics~\cite{CGS94}.  
As will be evident from \cref{thm:psabnd}, overestimating $d$ leads to 
lower (interior) bounds on $\PSA(\BA,\BE)$.)

Expand the solution of the DAE~\cref{eq:EmuEmu} in the Schur basis as
\[ \Bx(t) = \BQ_\mu \By(t) + \wt{\BQ}_\mu \Bz(t).\]
Substitute this form for $\Bx(t)$ and the Schur decomposition 
into~\cref{eq:EmuEmu} to arrive at the apparently coupled equations
\begin{align}
\BG_\mu \By'(t) + \BD_{\kern-.5pt \mu} \Bz'(t) &= (\BI+\mu@\BG_\mu)\By(t) + \cool{\mu}@\BD_{\kern-.5pt \mu} \Bz(t) \label{eq:dae2a}\\[0.25em]
              \BN_\mu \Bz'(t) &= (\BI+\mu@\BN_\mu)\Bz(t)\label{eq:dae2b}.
\end{align}
Since $\BN_\mu$ is nilpotent, $\BN_\mu^d = \Bzero$, 
premultiplying~\cref{eq:dae2b} by $\BN_\mu^{d-1}$ implies that, for all $t$,
\[\Bzero = \BN_\mu^d \Bz'(t) = (\BN_\mu^{d-1} + \mu@\BN_\mu^d) \Bz(t)
                             = \BN_\mu^{d-1} \Bz(t),\]
and hence $\Bzero = \BN_\mu^{d-1} \Bz'(t)$.
Thus $\Bzero = \BN_\mu^{d-1} \Bz'(t) = (\BN_\mu^{d-2} + \mu \BN_\mu^{d-1})\Bz(t) = \BN_\mu^{d-2} \Bz(t)$,
and so $\Bzero = \BN_\mu^{d-2}@\Bz'(t)$.
Repeating this process with decreasing powers of $\BN_\mu$
eventually leads to $\Bzero = \Bz(t)$
for all $t$.
Consequently equation~\cref{eq:dae2a} becomes
\[ \BG_\mu\By'(t) = (\BI+\mu@\BG_\mu)@\By(t).\]
Inverting $\BG_\mu$, we arrive at the solution
\begin{equation} \label{eq:dae3}
     \Bx(t) = \BQ_\mu^{}@\eop^{t(\BG_\mu^{-1}+\mu\BI)} \BQ_\mu^*@\Bx(0),
\end{equation}
with the stipulation that $\Bx(0) \in \Ran(\BQ_\mu)$
\new{to ensure that the initial condition is consistent with} the algebraic
constraints implicit in the DAE.\ \ 

Since $\|\Bx(t)\| = \|\eop^{t@(\BG_\mu^{-1}+\mu\BI)} (\BQ_\mu^*\Bx(0))\|$,
the solution~\cref{eq:dae3} suggests a definition for the 
$\eps$-pseudospectrum of the pencil~$(\BA,\BE)$ that is appropriate
for analyzing the transient behavior of DAEs, a direct generalization
of the approach commonly used for standard dynamical
systems.
We propose to define
\[ \PSA(\BA,\BE) 
   := \{z\in\C: \|(z\BI-(\BG_\mu^{-1}+\mu\BI))^{-1}\| > 1/\eps\}.\]
It appears this $\PSA(\BA,\BE)$ depends on~$\mu$, 
but since $\mu$ was just a device introduced to arrive at 
a solution formula, it should have no influence on the dynamics.  
Does $\mu$ affect these pseudospectra?

\subsection{Independence from \boldmath $\mu$}
Suppose for $\mu, \nu \in \C$ both $\BA-\mu\BE$ and $\BA-\nu\BE$ 
are invertible.  
The spectra of 
$\BE_\mu := (\BA-\mu@\BE)^{-1}\BE$ 
and
$\BE_\nu := (\BA-\nu\BE)^{-1}\BE$ 
are closely related.
Suppose $\lambda \in \SPEC(\BE_\mu)$, so for some nonzero $\Bx\in\Cn$,
$(\BA-\mu@\BE)^{-1} \BE@\Bx = \lambda \Bx$.
Thus
\begin{align*}
\BE@\Bx &= \lambda (\BA-\mu@\BE) \Bx \\[.25em]
        &= \lambda (\BA-\nu@\BE)(\BI + (\nu-\mu)(\BA-\nu\BE)^{-1}\BE)\Bx.
\end{align*} 
Premultiply by $(\BA-\nu@\BE)^{-1}$ to get
\[ \BE_\nu \Bx = \lambda(\BI + (\nu-\mu)@\BE_\nu)@\Bx,\]
which is equivalent to
\[ (1 + \lambda(\mu-\nu))@\BE_\nu\Bx = \lambda \Bx.\]
Notice that $1+\lambda(\mu-\nu)= 0$ would imply both $\lambda\ne 0$ and
$\Bzero = \lambda \Bx$; since $\Bx\ne \Bzero$, this is impossible.
Thus we have
\[ \BE_\nu\Bx = {\lambda \over 1 + \lambda(\mu-\nu)} \Bx,\]
proving the following lemma.

\smallskip
\begin{lemma}
Suppose for $\mu, \nu \in \C$ both $\BA-\mu@\BE$ and $\BA-\nu\BE$ 
are invertible.  If $\lambda\in\SPEC(\BE_\mu)$, then
\[ {\lambda \over 1+\lambda(\mu-\nu)} \in \SPEC(\BE_\nu).\]
\end{lemma}
\smallskip

From the fact that $\BA-\mu\BE = \BA-\nu\BE + (\nu-\mu)\BE$ follows
\begin{equation} \label{eq:resid}
 (\BA-\mu\BE)^{-1} = (\BA-\nu\BE)^{-1} + (\mu-\nu)(\BA-\nu\BE)^{-1}\BE(\BA-\mu\BE)^{-1},
\end{equation}
a generalization of the ``first resolvent identity'' in standard spectral theory.
The Schur factorization $\Emu = (\BA-\mu\BE)^{-1}\BE = \BQ\BT\BQ^*$ and 
the identity~\cref{eq:resid} give 
\begin{align*}
  \BQ\BT &= (\BA-\mu@\BE)^{-1}\BE\BQ \\
         &= (\BA-\nu@\BE)^{-1}\BE\BQ + (\mu-\nu)(\BA-\nu\BE)^{-1}\BE(\BA-\mu\BE)^{-1}\BE\BQ.
\end{align*}
Substituting $(\BA-\mu\BE)^{-1}\BE = \BQ\BT\BQ^*$ on the right-hand side then gives
\[ \BQ\BT = (\BA-\nu\BE)^{-1}\BE\BQ(\BI+(\mu-\nu)\BT).\]
Since $1+(\mu-\nu)\lambda \ne 0$ for all eigenvalues $\lambda$ of $\Emu$,
$\BI+(\mu-\nu)\BT$ is invertible, so
\[ (\BA-\nu\BE)^{-1}\BE\BQ = \BQ\BT@(\BI+(\mu-\nu)\BT)^{-1}.\]
Note that $\BT@(\BI+(\mu-\nu)\BT)^{-1}$, the product of triangular
matrices, must itself be triangular, and 
$\BE_\nu := (\BA-\nu\BE)^{-1}\BE$ has the same Schur basis $\BQ$ 
as $\BE_\mu$.
Partition~$\BQ$ and $\BT$ as in~\cref{eq:schur}, so that
\[ \BQ^*\BE_\nu\BQ = \BT@(\BI+(\mu-\nu)\BT)^{-1}
   =    \left[\begin{array}{cc} \BG_\mu & \BD_{\kern-.5pt \mu} \\ \Bzero & \BN_\mu \end{array}\right]
        \left[\begin{array}{cc} \BI+(\mu-\nu)\BG_\mu & (\mu-\nu)\BD_{\kern-.5pt \mu} \\
                 \Bzero & \BI+(\mu-\nu)\BN_\mu \end{array}\right]^{-1}
\]
has $(1,1)$ block equal to 
\begin{equation} \label{eq:Gnu}
 \BG_\nu := \BG_\mu(\BI+(\mu-\nu)\BG_\mu)^{-1}
\end{equation}
and $(2,2)$ block equal to
\[ \BN_\nu := \BN_\mu(\BI+(\mu-\nu)\BN_\mu)^{-1}.\]
Since the eigenvalues of $\BG_\mu$ are nonzero, so too are those
of $\BG_\nu$.
Notice that $\BN_\nu^d = \BN_\mu^d (\BI+(\mu-\nu)\BN_\mu)^{-d}$
(as a function of $\BN_\mu$ commutes with $\BN_\mu$), 
so $\BN_\nu$ is also nilpotent.

Inverting both sides of~\cref{eq:Gnu} gives
$\BG_\nu^{-1} = (\BI + (\mu-\nu)\BG_\mu^{})\BG_\mu^{-1}$, 
which simplifies to
\begin{equation} \label{eq:indpt}
 \BG_\nu^{-1} + \nu \BI = \BG_\mu^{-1} + \mu@\BI.
\end{equation}
It follows that $\BG_\mu^{-1} + \mu@\BI$ is independent of $\mu$ 
(provided that $\BA-\mu@\BE$ is invertible),
allowing us to sharpen up our definition of the 
$\eps$-pseudospectrum of a matrix pencil.

\smallskip
\begin{definition} \label{def:psa}
Suppose $(\BA,\BE)$ is a regular matrix pencil, and $\mu\in\C$ is 
any value for which $\BA-\mu@\BE$ is invertible.  
Let $\BG_\mu$ be the submatrix in the Schur factorization~$\cref{eq:schur}$ 
corresponding to the nonzero eigenvalues of $(\BA-\mu@\BE)^{-1}\BE$,
\cool{and $\|\cdot\|$ denote the 2-norm.}
For any $\eps>0$, the $\eps$-pseudospectrum of $(\BA,\BE)$ is defined to be
\begin{align}
 \PSA(\BA,\BE) &:= \{z\in\C: \|(z\BI-(\BG_\mu^{-1}+\mu@\BI))^{-1}\| > 1/\eps\} 
   \label{eq:daepsa1} \\[.25em]
               &\phantom{:}= \{z\in\C: \|((z-\mu)\BG_\mu - \BI)^{-1}\BG_\mu^{}\| > 1/\eps\} 
   \nonumber \\[.25em]
               &\phantom{:}= \PSA(\BG_\mu^{-1}) + \mu,
   \nonumber
\end{align}
where $\PSA(\BG_\mu^{-1})$ refers to
the standard matrix $\eps$-pseudospectrum~$\cref{eq:psa1}$ of $\BG_\mu^{-1}$.
The set $\PSA(\BA,\BE)$ is independent of $\mu$.
\end{definition}
\smallskip

\begin{figure}[b!]
\begin{center}
\includegraphics[scale=0.4]{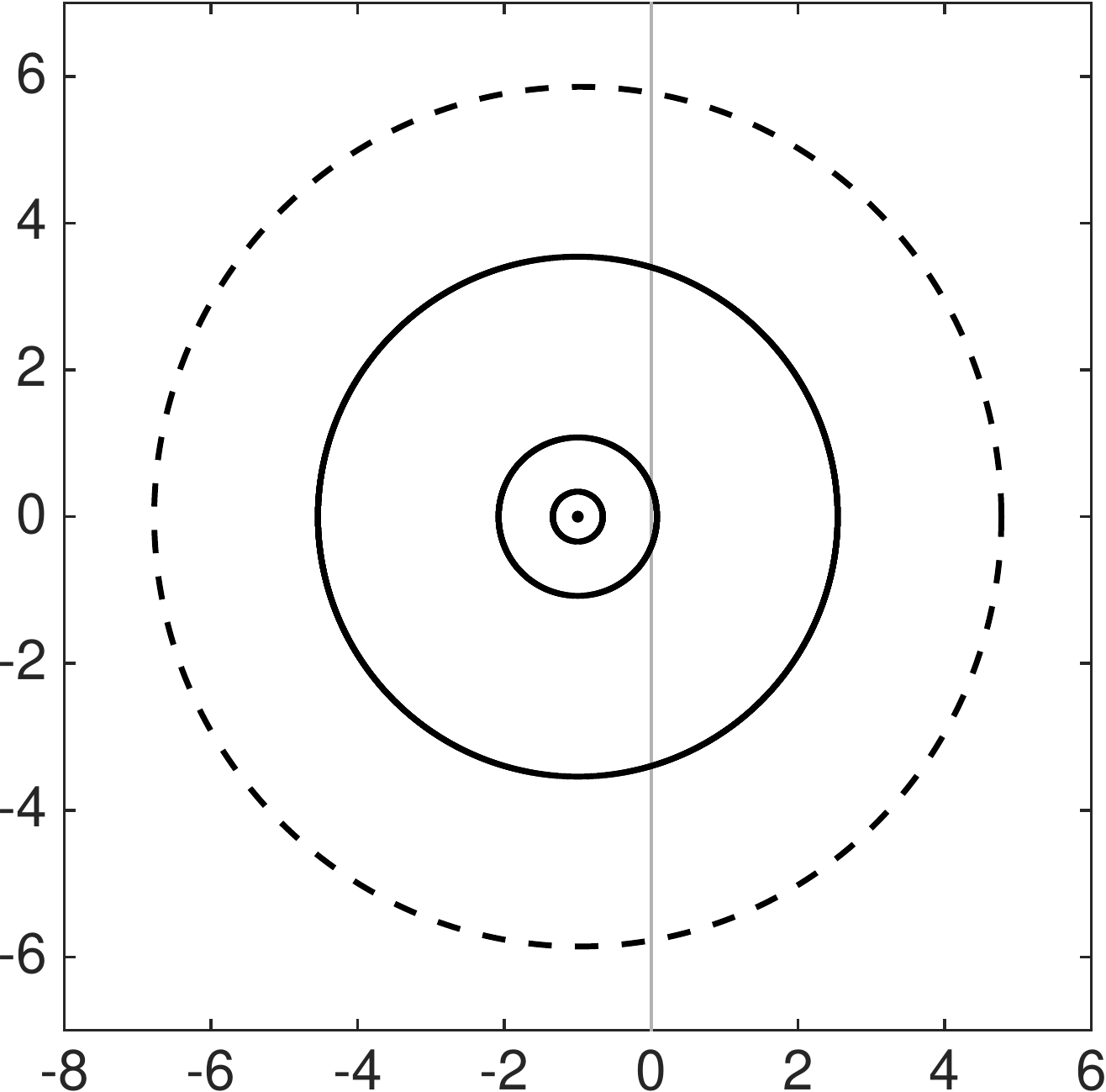}\qquad\qquad
\includegraphics[scale=0.4]{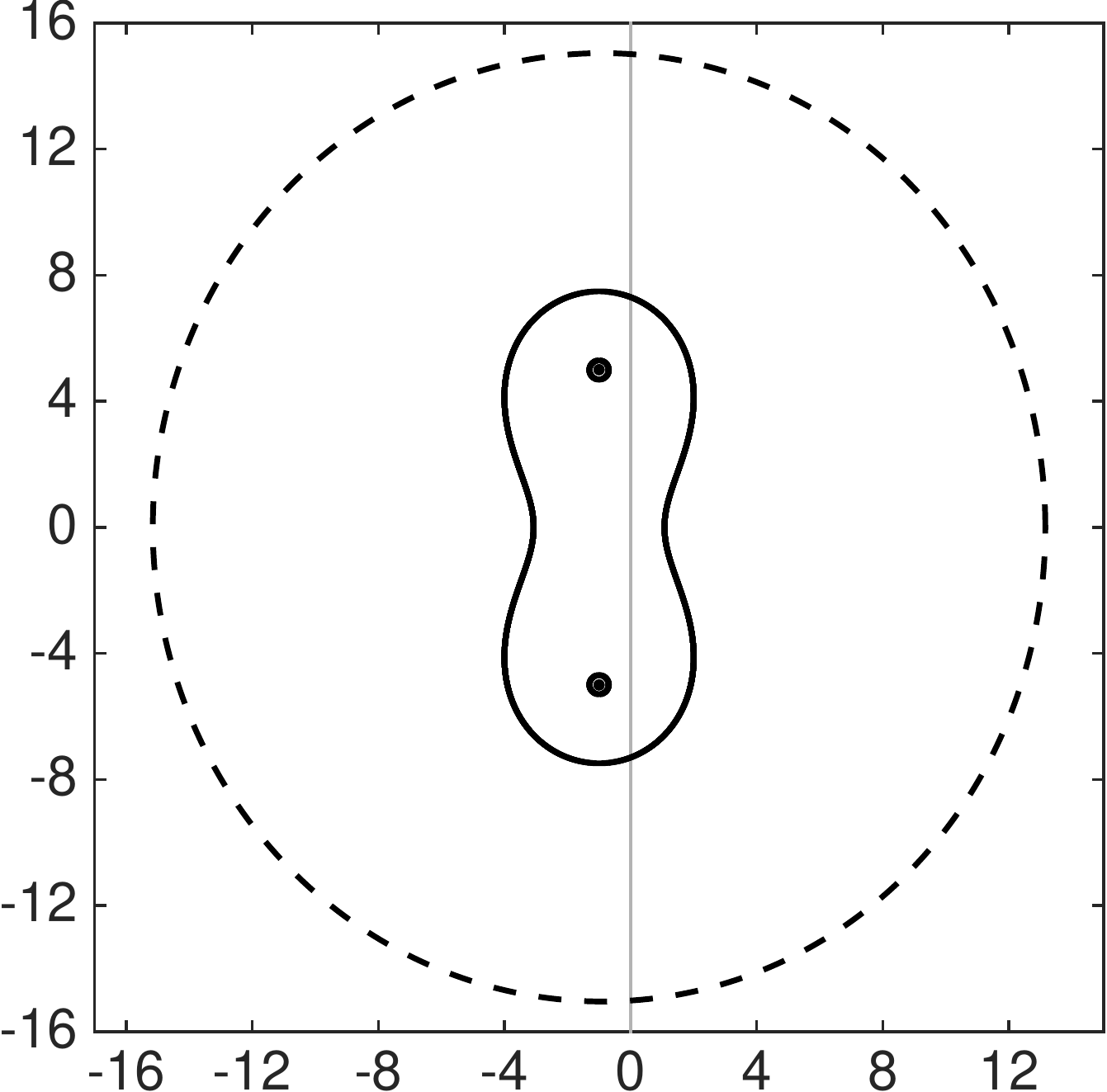}
\end{center}
\caption{ \label{fig:psa1}
Boundaries of $\eps$-pseudospectra for $(\BA,\BE)$ using \cref{def:psa}
for $\eps = 10^{0}$, $10^{-1}$, and $10^{-2}$ (solid curves) and $W(\BA,\BE)$ 
using \cref{def:nr} (dashed curves); $(\BA,\BE)$ are the same as for
the left and right plots in \cref{fig:ex3d}.
}
\end{figure}

\Cref{fig:psa1} shows pseudospectra, as defined by \cref{def:psa},
for the pairs~$(\BA,\BE)$ used in \cref{fig:ex3d}.  
In both cases the $\eps=1$ pseudospectrum contains points 
$z$ for which $\Re\,z > \eps$, which, as we shall see in the next section,
guarantees the solution $\Bx(t)$ to the DAE~\cref{eq:dae}
exhibits transient growth for some (valid) initial condition $\Bx(0)$.
More sophisticated examples appear in \cref{sec:examples}.

\begin{remark}   \rm
We collect several observations about this definition.
\begin{remunerate}
\item In the Schur decomposition~\cref{eq:schur}, 
the nonzero eigenvalues can be rearranged in any order on the 
diagonal of $\BG_\mu$; this reordering effectively replaces $\BG_\mu$ 
with some unitary similarity transformation, $\BU^*\BG_\mu\BU$.  
By the unitary invariance of the 2-norm, this transformation will not 
affect the definition of $\PSA(\BA,\BE)$.  
Beware, though, that if one 
independently computes Schur decompositions of $\BE_\mu$ and $\BE_\nu$ 
for $\mu\ne \nu$, one will likely find 
$\BG_\mu^{-1}+\mu\BI \ne \BG_\nu^{-1} + \nu\BI$ 
due to such a unitary similarity transformation.
\item 
When $\BE$ is invertible, \cref{def:psa} reduces to 
Ruhe's definition~\cref{eq:ruhe}.
(Take $\mu=0$ in the definition, and again use unitary invariance
of the 2-norm.)
\item Practically speaking, $\mu$ should be chosen so that 
$\BE_\mu = (\BA-\mu\BE)^{-1}\BE$ (and its Schur factor, from 
which we extract $\BG_\mu$) can be computed reliably.  
In cases where $\BA$ is invertible and well conditioned,
$\mu=0$ is a natural choice.  
For large-scale problems like those in \cref{sec:examples},
$\mu$ should be chosen to influence the convergence of a projection
method for computing an approximate invariant subspace.
\item Since $\PSA(\BA,\BE) = \PSA(\BG_\mu^{-1}+\mu@\BI)$ 
is just a standard pseudospectrum, one can compute these sets using the
algorithms and software packages designed for standard pseudospectra;
see, e.g., \cite[chaps.~39--42]{TE05}, \cite{Wri02b,Wri02a}.
\item Note that the angle between the invariant subspaces associated 
with the finite and infinite eigenvalues (controlled by the off-diagonal
block $\BD$) does not influence this definition of pseudospectra,
just as it does not play a role in the solution $\Bx(t)$ of the
DAE in~\cref{eq:dae3}.
Were $\BE$ perturbed slightly to become invertible (say, $\BE \to \BE+\delta \BI$),
$\BD$ would certainly influence the pseudospectra $\PSA(\BE^{-1}\BA)$, 
just as such a perturbation would remove the algebraic constraint
on $\Bx(0)$ and allow initial conditions with components in the 
invariant subspace previously associated with the infinite eigenvalues.
\item Similarly, the index of the DAE (revealed through the degree of
nilpotency of $\BN_\mu$) influences neither $\Bx(t)$ nor $\PSA(\BA,\BE)$.

\end{remunerate}
\end{remark}

\subsection{Numerical range} \label{sec:nr}
We can similarly generalize the definition of the
numerical range (field of values) of a matrix $\BA\in\Cnn$,
\[ \NR(\BA) := \left\{ \Bx^*\BA\Bx: \Bx\in\C^n, \|\Bx\|=1\right\}.\]

\begin{definition} \label{def:nr}
The \emph{numerical range} (or \emph{field of values}) of the regular matrix pencil $(\BA,\BE)$,
\cool{in the Euclidean inner product and the 2-norm},
is
\begin{align*}
    \NR(\BA,\BE) &:= \left\{ \By^*\BG_\mu^{-1}\By + \mu: \By\in\C^{n-d}, \|\By\|=1\right\} \\[.25em]
                 &\phantom{:}= \NR(\BG_\mu^{-1}) + \mu,
\end{align*}
where $\mu\in\C$ is any value for which $\BA-\mu@\BE$ is invertible.  
\end{definition}

\medskip
\Cref{fig:psa1} shows $W(\BA,\BE)$ for the same 
matrices used in the earlier examples.
Like our proposal for pseudospectra,
this definition for the numerical range 
differs from the conventional approach 
for matrix pencils~\cite{Hoc11,LR94,Psa00},
but, as we will see in \cref{sec:transient}, 
it gives important insight into transient dynamics.
With our definitions, $\PSA(\BA,\BE)$ can be bounded
in terms of $\NR(\BA,\BE)$.

\begin{theorem}
Let $(\BA,\BE)$ be a regular pencil.  For all $\eps>0$, 
\begin{equation} \label{eq:stone}
  \PSA(\BA,\BE) \subseteq \NR(\BA,\BE) + \{z\in\C: |z| < \eps\}.
\end{equation}
\end{theorem}
\begin{proof}
Let $\mu\in\C$ be any value for which $\BA-\mu@\BE$ is invertible.
Then $\PSA(\BA,\BE) = \PSA(\BG_\mu^{-1}) + \mu$ and 
$\NR(\BA,\BE) = \NR(\BG_\mu^{-1})+\mu$.
The inclusion~\cref{eq:stone} then follows by applying the analogous
bound for matrices:  $\PSA(\BG_\mu^{-1})\subseteq \NR(\BG_\mu^{-1})+\{z\in\C:|z|<\eps\}$; see \cite[thm.~4.20]{Sto32}, \cite[p.~169]{TE05}.\hfill
\end{proof}

\subsection{Other norms} \label{sec:norms}
In many applications, one seeks to measure the transient behavior of $\Bx(t)$
not in the vector 2-norm, but in some norm that has more physical relevance.
For example, $\|\Bx(t)\|^2$ could measure the instantaneous energy in a system.
When the system is a discretized partial differential equation,
the norms should not bear any $n$-dependence as the discretization is refined.
\Cref{sec:examples} gives a specific example from fluid dynamics,
where subvectors of $\Bx(t)$ must be handled differently.

For clarity, in this subsection we use the notation $\|\cdot\|_2$ and $\sigma_{\eps,2}(\cdot)$,  
while the subscript~``2'' is implicit in the notation $\|\cdot\|$ and $\sigma_\eps(\cdot)$ \cool{elsewhere in this section}.
Let $\BH\in\Cnn$ be a Hermitian positive definite matrix factored
as $\BH = \BR^*\BR$ for some $\BR\in\Cnn$
(e.g., $\BR$ is a Cholesky factor or the Hermitian square root of $\BH$).
Consider the inner product $\langle \cdot, \cdot \rangle_\BH$ 
defined for $\Bx,\By\in\Cn$ by
\[ \langle \Bx, \By\rangle_\BH := \By^*\BH\Bx = (\BR\By)^*(\BR\Bx).\]
This inner product induces the vector norm 
\begin{equation} \label{eq:Hvecnorm}
 \|\Bx\|_\BH := \langle \Bx,\Bx\rangle_\BH^{1/2} = \|\BR\Bx\|_2,
\end{equation}
with which we associate, for any $\BM\in\Cnn$, the matrix norm
\begin{align} 
   \|\BM\|_\BH &:= \max_{\Bx\ne \Bzero} {\|\BM\Bx\|_\BH \over \|\Bx\|_\BH}
                    \nonumber  \\[3pt] 
               &= \max_{\Bx\ne \Bzero} {\|\BR\BM\Bx\|_2 \over \|\BR\Bx\|_2} 
                = \max_{\Bx\ne \Bzero} {\|\BR\BM\BR^{-1}(\BR\Bx)\|_2 \over \|\BR\Bx\|_2} 
                = \|\BR\BM\BR^{-1}\|_2. \label{eq:Hmatnorm}
\end{align}
The $\BH$-norm of $\BM$ is just the 2-norm of the similar matrix $\BR\BM\BR^{-1}$,
giving a simple way to compute $\|\BM\|_\BH$.

The definition of the standard matrix pseudospectrum $\sigma_\eps(\BA)$ easily
accommodates any norm induced by a general inner product:
simply use $\|\cdot\|_\BH$ for the norm in (\ref{eq:psa1})--(\ref{eq:psa2}).
Via the calculation~(\ref{eq:Hmatnorm}), one can use software for
2-norm pseudospectra (e.g., EigTool~\cite{Wri02a})
to compute $\BH$-norm pseudospectra, since
$\sigma_{\eps,\BH}(\BA) = \sigma_{\eps,2}(\BR\BA\BR^{-1})$.
Adapting Definition~\ref{def:psa} for the matrix pencil pseudospectrum
$\sigma_\eps(\BA,\BE)$ to 
a norm $\|\cdot\|_\BH$ induced by a general inner product requires more care.%
\footnote{Theoretically the matter is trivial: require
$[\BQ_\mu\ \wt{\BQ}_\mu]$ in the Schur decomposition~(\ref{eq:schur})
to be unitary with respect to the $\BH$-inner product,
and replace $[\BQ_\mu\ \wt{\BQ}_\mu]^*$ in the analysis with the
$\BH$-adjoint of $[\BQ_\mu\ \wt{\BQ}_\mu]$.
We provide a more concrete discussion for computational convenience.}
We discuss two equivalent approaches.

\subsubsection{Approach~1: Transform state vector coordinates}
We seek to measure transient behavior of the DAE
solution $\Bx(t)$ in the $\BH$-norm.
By~\cref{eq:Hvecnorm}, $\|\Bx(t)\|_\BH = \|\BR\Bx(t)\|_2$.
Substituting $\Bs(t) := \BR\Bx(t)$ into~\cref{eq:dae} 
leads to the DAE
\[ \BE\BR^{-1}\Bs'(t) = \BA\BR^{-1}\Bs(t),\]
suggesting that one simply define
\begin{equation} \label{eq:psaH}
 \sigma_{\eps,\BH}(\BA,\BE) := \sigma_{\eps,2}(\BA\BR^{-1},\BE\BR^{-1}).
\end{equation}
This definition behaves as expected when $\BE$ is invertible:
$\sigma_{\eps,\BH}(\BA,\BE)$, as given in~\cref{eq:psaH}, 
reduces to the $\BH$-norm pseudospectrum of $\BE^{-1}\BA$:
\begin{align*}\sigma_{\eps,\BH}(\BA,\BE) 
    &= \sigma_{\eps,2}(\BA\BR^{-1}, \BE\BR^{-1})  \\
    &= \sigma_{\eps,2}((\BE\BR^{-1})^{-1}(\BA\BR^{-1})) 
    = \sigma_{\eps,2}(\BR\BE^{-1}\BA\BR^{-1})
    = \sigma_{\eps,\BH}(\BE^{-1}\BA).
\end{align*}
For singular $\BE$, definition~\cref{eq:psaH}
involves a Schur factorization of 
\begin{equation} \label{eq:EmuHdef}
 \BE_{\mu,\BH} := (\BA\BR^{-1}-\mu@\BE\BR^{-1})^{-1}\BE\BR^{-1} 
                  = \BR(\BA-\mu@\BE)^{-1}\BE\BR^{-1},
\end{equation}
which can be partitioned in the form~\cref{eq:schur}.
The $(1,1)$ block of the central factor in this decomposition,
denoted~$\BG_\mu$ in~\cref{eq:schur}, generally depends on $\BR$.

\subsubsection{Approach~2: Transform the Schur factorization~\cref{eq:schur}}
\label{sec:Hnorm2}
Suppose one has a Schur factorization~\cref{eq:schur}
in the Euclidean inner product for $\BE_\mu = (\BA-\mu@\BE)^{-1}\BE$.
How does $\BG_\mu$, key to \cref{def:psa}, 
change with the inner product?
Using~\cref{eq:schur},
\begin{equation} \label{eq:EmuH}
 \BE_{\mu,\BH} = \BR \BE_\mu \BR^{-1} 
       =  \BR 
          \left[\begin{array}{cc} \BQ_\mu & \wt{\BQ}_\mu\end{array}\right]
          \left[\begin{array}{cc} \BG_\mu & \BD_{\kern-.5pt \mu} \\ \Bzero & \BN_\mu\end{array}\right]
          \left[\begin{array}{c} \BQ_\mu^* \\ \wt{\BQ}_\mu^*\end{array}\right] \BR^{-1}.
\end{equation}
Compute a QR factorization 
\[ \BR \left[\begin{array}{cc}\BQ_\mu & \wt{\BQ}_\mu \end{array}\right]
       = \left[\begin{array}{cc}\BZ_\mu & \wt{\BZ}_\mu \end{array}\right]
          \left[\begin{array}{cc}\BS_\mu & \times \\ \Bzero & \wt{\BS}_\mu \end{array}\right],
\]
where the first matrix on the right is unitary, and $\times$ is a generic
placeholder for a submatrix that does not factor into our discussion. 
Note that the columns of $\BZ_\mu$ form an orthonormal basis for the range of $\BR\BQ_\mu$.
Substituting the QR factorization into~\cref{eq:EmuH} gives
\begin{equation} \label{eq:schurH}
 \BE_{\mu,\BH} 
       =  \left[\begin{array}{cc} \BZ_\mu & \wt{\BZ}_\mu\end{array}\right]
          \left[\begin{array}{cc} \BS_{\mu}^{} \BG_\mu^{} \BS_{\mu}^{-1} & \times \\
                   \Bzero & \wt{\BS}_\mu^{} \BN_\mu^{} \wt{\BS}_\mu^{-1}\end{array}\right]
          \left[\begin{array}{cc} \BZ_\mu^* \\ \wt{\BZ}_\mu^*\end{array}\right].
\end{equation}
This analogue of~\cref{eq:schur} reveals how the $\BH$-inner product 
affects the pseudospectra:

\begin{center}
\vspace*{4pt}
\begin{tabular}{rl}
2-norm $\eps$-pseudospectrum:  
           & $\sigma_{\eps,2}(\BA,\BE) = \sigma_{\eps,2}(\BG_\mu^{-1})+\mu$ \\[5pt]
$\BH$-norm $\eps$-pseudospectrum: 
           & $\sigma_{\eps,\BH}(\BA,\BE) = \sigma_{\eps,2}(\BS_\mu^{}\BG_\mu^{-1}\BS_\mu^{-1})+\mu$.
\end{tabular}
\vspace*{4pt}
\end{center}
The situation perfectly parallels the case of invertible $\BE$: 
in that case, the $\BH$-norm pseudospectra of $\BE^{-1}\BA$ 
are the 2-norm pseudospectra of a similarity transformation with $\BR$.
For singular $\BE$, this similarity transformation is not with $\BR$,
but with $\BR$ filtered through the subspace $\Ran(\BQ_\mu)$ in which
the solution evolves.

In summary, to compute $\sigma_{\eps,\BH}(\BA,\BE)$:\\[-8pt]
\begin{enumerate}
\item Compute the Schur factorization~\cref{eq:schur} 
of $\BE_\mu := (\BA-\mu@\BE)^{-1}\BE$ to get $\BG_\mu$, $\BQ_\mu$.
\item Compute the economy-sized QR factorization $\BR\BQ_\mu = \BZ_\mu\BS_\mu$.
\item Compute $\sigma_{\eps,\BH}(\BA,\BE) = \sigma_{\eps,2}(\BS_\mu^{}\BG_\mu^{-1}\BS_\mu^{-1})+\mu$.
\end{enumerate}

\subsubsection{Norms not induced by inner products}

We shall not dwell long on norms that are not induced by inner products. 
The solution formula~\cref{eq:dae3} still holds, 
so $\|\Bx(t)\| = \|\BQ_\mu^{}@\eop^{t(\BG_\mu^{-1}+\mu\BI)} \BQ_\mu^*@\Bx(0)\|$.
Given a system of submultiplicative norms, 
\begin{equation} \label{eq:banach}
 \|\Bx(t)\| \le \|\BQ_\mu\|@\|\BQ_\mu^*\| \|\eop^{t(\BG_\mu^{-1}+\mu\BI)}\| \|\Bx_0\|.
\end{equation}
For example, for the matrix 1-norm, $\|\BQ_\mu\| \le \sqrt{n}$ and 
$\|\BQ_\mu^*\| \le \sqrt{n-d}$.
Thus \cref{def:psa} can still be justified (for example, 
$\BD_{\kern-.5pt \mu}$ in~\cref{eq:schur} plays no role in $\|\Bx(t)\|$, and so
should not factor in $\PSA(\BA,\BE)$),  
but the additional constants in~\cref{eq:banach}
make the resulting bounds less satisfying
than those for norms induced by inner products.

\section{Transient behavior} \label{sec:transient}
Throughout this section we assume that $(\BA,\BE)$ is 
asymptotically stable, i.e., all finite eigenvalues of 
the pencil fall strictly in the left-half plane, and
hence $\Bx(t) \to \Bzero$ as $t\to\infty$ for all 
$\Bx(0)$ that satisfy the algebraic constraints imposed
by the DAE.\ \ 
\cool{We use the 2-norm here, but \cref{sec:norms} makes clear
how the results that follow can be adapted to any norm defined
by an inner product.}
We seek to identify situations where
$\|\Bx(t)\|$ grows before its asymptotic decay 
(or converges more slowly than would be predicted from the
pencil's rightmost finite eigenvalue), 
as shown in \cref{fig:ex3d}.
\Cref{def:psa,def:nr} were designed to
illuminate this transient behavior.

As usual, let $\mu\in\C$ be such that $\BA-\mu@\BE$ is invertible.
Using the notation of the last section,
any valid initial condition for the DAE must satisfy $\Bx(0) \in \Ran(\BQ_\mu)$,
and hence can be written as $\Bx(0) = \BQ_\mu \By_0$ for some $\By_0 \in \C^{n-d}$.
Using the unitary invariance of the norm, 
\begin{align*}
    \norm{\Bx(t)}  &= \norm{\BQ_\mu^{} \eop^{t(\BG_\mu^{-1}+\mu\BI)}\BQ_\mu^* \Bx(0)} \\[0.25em]
                   &= \norm{ \eop^{t(\BG_\mu^{-1}+\mu\BI)} \By_0}.
\end{align*}
Similarly, since for any $\By_0\in\C^{n-d}$, 
$\BQ_\mu\By_0$ is a valid initial condition for the DAE,
the definition of the matrix norm implies that for any $t$,
there exists some unit vector $\Bx(0)\in\Ran(\BQ_\mu)$ such that
\[ \norm{\Bx(t)} = \norm{\eop^{t(\BG_\mu^{-1}+\mu\BI)}}.\]
We thus have available the wealth of results characterizing
the transient behavior of $\Bx(t)$ based on spectral properties 
of $\BG_\mu^{-1}+\mu@\BI$.
We state a number of bounds that now
follow as easy corollaries of results for standard dynamical systems.
For conventional pseudospectra, proofs of these results can be 
found in~\cite[part~IV]{TE05}.
We first define the key quantities that connect pseudospectra and the 
numerical range to the
transient behavior of continuous-time systems.

\begin{definition} \label{def:psabs}
The \emph{$\eps$-pseudospectral abscissa} of the regular pencil $(\BA,\BE)$ 
is 
\[ \psabs(\BA,\BE) := \sup_{z\in \PSA(\BA,\BE)} \Re\,z.\]
\end{definition}

\begin{definition} \label{def:nabs}
The \emph{numerical abscissa} of the regular pencil $(\BA,\BE)$ 
is 
\[ \nabs(\BA,\BE) := \sup_{z\in \NR(\BA,\BE)} \Re\,z.\]
\end{definition}

\medskip
The analogue of $\nabs(\BA,\BE)$ in the standard matrix case is sometimes 
called the \emph{logarithmic norm}~\cite{Sod06}.
Note that $\psabs(\BA,\BE)$ and $\nabs(\BA,\BE)$ can be computed from their
standard matrix analogues: 
\begin{equation} \label{eq:psnabs}
 \psabs(\BA,\BE)  = \psabs(\BG_\mu^{-1}) + \mu, 
    \qquad 
   \nabs(\BA,\BE)  = \nabs(\BG_\mu^{-1}) + \mu,
\end{equation}
with both quantities independent of $\mu$.
The latter equality implies 
\[ \nabs(\BA,\BE) = \mu + \lambda_{\rm max}\Big({\BG_\mu^{-1} + \BG_\mu^{-*}\over 2}\Big),\] 
where $\lambda_{\rm max}(\cdot)$ is the rightmost eigenvalue of a Hermitian matrix
\cool{\cite[lemma~1.5.7]{HJ91}}.

\subsection{Behavior at \boldmath $t=0$}
The numerical range describes the early behavior of a 
dynamical system, limiting the rate at which $\|\Bx(t)\|$ can initially grow.

\begin{theorem} \label{thm:t0}
Let $(\BA,\BE)$ be a regular pencil with $\BA-\mu@\BE$ invertible.
Then 
\[ {\dop \over \dop t} \big\|\eop^{t(\BG_\mu^{-1}+\mu\BI)}\big\| \bigg|_{t=0} = \nabs(\BA,\BE).\]
For any unit vector $\Bx(0) \in \Ran(\BQ_\mu)$, the solution $\Bx(t)$ to 
$\BE@\Bx'(t) = \BA\Bx(t)$ thus satisfies
\[ {\dop \over \dop t} \|\Bx(t)\| \bigg|_{t=0} \le \nabs(\BA,\BE),\]
with equality attained for some unit vector $\Bx(0) \in \Ran(\BQ_\mu)$.
\end{theorem}

See~\cite[chap.~17]{TE05} for a proof in the standard matrix case,
which can be applied to $\BG_\mu^{-1}+\mu\BI$ to obtain \cref{thm:t0}.
This result is connected to the
Lumer--Phillips theorem, which relates dissipative operators to contraction 
semigroups~\cite[sect.~1.4]{Paz83}.  

If $\nabs(\BA,\BE)>0$, the system must exhibit transient growth for some initial conditions.
The maximum growth rate is attained for $\Bx(0) = \BQ_\mu \By$, where $\By$ is a
unit eigenvector associated with the rightmost eigenvalue of $\BG_\mu^{-1}+\BG_\mu^{-*}$.

\subsection{Lower bounds on maximal growth}

When $\nabs(\BA,\BE)>0$, the numerical range 
captures the initial growth of $\|\Bx(t)\|$, 
but it does not address the extent of that growth at times $t>0$.  
Pseudospectra are more useful for this task.
The next theorem implies that if $\PSA(\BA,\BE)$ extends
more than $\eps$ into the right-half plane, then there exists some $\Bx(0)$ for 
which $\Bx(t)$ grows by at least a factor of $\psabs(\BA,\BE)/\eps$.

\begin{theorem} \label{thm:explow}
Let $(\BA,\BE)$ be a regular pencil with $\BA-\mu@\BE$ invertible.  
Then 
\begin{equation} \label{eq:laplow}
 \sup_{t\ge 0} \big\|\eop^{t(\BG_\mu^{-1}+\mu\BI)}\big\| 
    \ge {\alpha_\eps(\BA,\BE) \over \eps}
\end{equation}
for all $\eps>0$,
and there exists some $\Bx(0) \in \Ran(\BQ_\mu)$ such that the solution 
$\Bx(t)$ to $\BE@\Bx'(t) = \BA\Bx(t)$
realizes this transient growth:
\[ \sup_{t\ge 0}\ {\|\Bx(t)\| \over \|\Bx(0)\|} \ge {\alpha_\eps(\BA,\BE)\over \eps}.\]
\end{theorem}

\begin{figure}[b!]
\begin{center}
\includegraphics[scale=0.45]{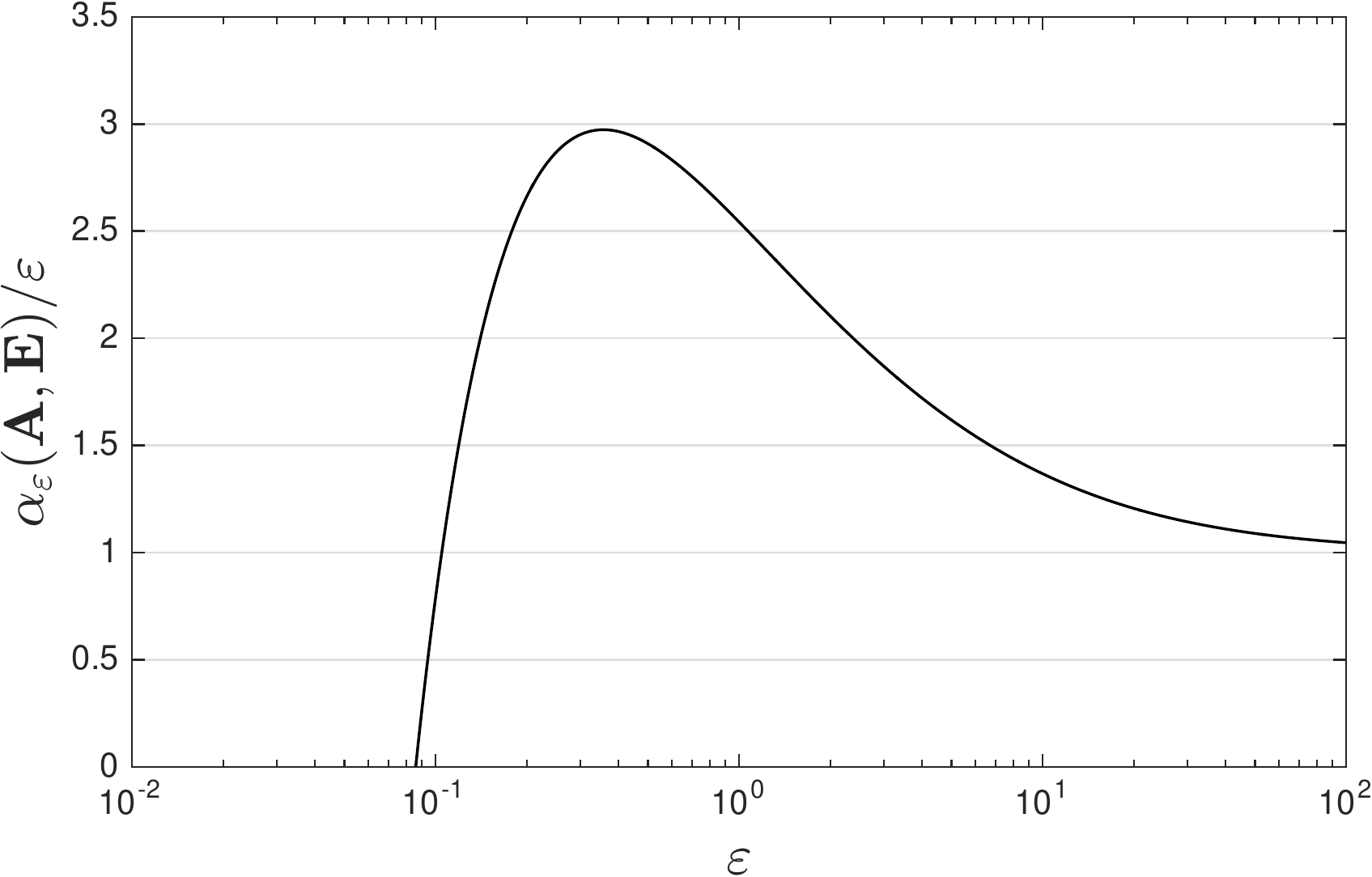}
\end{center}

\vspace*{-8pt}
\caption{\label{fig:ex1sabs}
The ratio $\alpha_\eps(\BA,\BE)/\eps$ as a function of $\eps$ for the 
example in~$\cref{eq:ex3da}$ for which $\PSA(\BA,\BE)$ was 
plotted on the left side of \cref{fig:psa1}.  
By \cref{thm:explow} there exists an initial condition
$\Bx(0) \in \Ran(\BQ_0)$ such that $\|\Bx(t)\|$ 
grows at least by a factor of nearly~3 (since ${\cal K}(\BA,\BE)$, 
the maximum of $\psabs(\BA,\BE)/\eps$, is nearly~3).
}
\end{figure}

The proof is a simple consequence of the identity equating the resolvent 
to the Laplace transform of the exponential of a matrix; see, e.g.,~\cite[thm.~11$\eps$]{ET01}.
\Cref{fig:ex1sabs} shows $\psabs(\BA,\BE)/\eps$ as a function of $\eps$
for the pencil in~\cref{eq:ex3da} whose pseudospectra were 
shown in the left plot of~\cref{fig:psa1}.

The $\eps$ that gives the greatest lower bound in \cref{thm:explow}
is of special interest.

\begin{definition} \label{def:kreiss}
The \emph{Kreiss constant} (with respect to the left-half plane)
of the regular pencil $(\BA,\BE)$ is 
\[ {\cal K}(\BA,\BE) :=  \sup_{\eps>0} {\alpha_\eps(\BA,\BE) \over \eps}.\]
\end{definition}

\Cref{thm:explow} is the most useful lower bound on 
transient growth, but it does not mark the time at which that growth
is realized.  
Some sense of time scale follows by adapting a bound of Trefethen
for the standard case~\cite[eq.~(14.13)]{TE05}.

\begin{theorem} \label{thm:lnt}
Let $(\BA,\BE)$ be a regular pencil with $\BA-\mu@\BE$ invertible,
and suppose that $\alpha_\eps(\BA,\BE)>0$ for some given $\eps>0$.  
Then for all $\tau>0$,
\begin{equation} \label{eq:lnt}
 \max_{t\in [0,\tau]} \|\eop^{t(\BG_\mu^{-1}+\mu\BI)}\| 
   \,\ge\, \eop^{\tau \psabs(\BA,\BE)} 
        \bigg( {1 \over 
                1 + \eps\big(\eop^{\tau \psabs(\BA,\BE)}-1\big)/\psabs(\BA,\BE)}\bigg),
\end{equation}
and for \cool{each $\tau>0$ there exists some initial condition $\Bx(0)\in\Ran(\BQ_\mu)$ and $t\in[0,\tau]$ such that}  $\|\Bx(t)\|/\|\Bx(0)\|$ attains this growth.
\end{theorem}

\subsection{Upper bounds on transient growth}
We now turn to upper bounds on $\|\Bx(t)\|$.  
The simplest bound, sometimes called Coppell's inequality
in the standard matrix case~\cite[sect.~4.2.1]{GQ95},
uses the numerical abscissa to limit the extent of growth at any given $t\ge0$.  

\begin{theorem}
Let $(\BA,\BE)$ be a regular pencil with $\BA-\mu\BE$ invertible.  Then
\[\big\|\eop^{t(\BG_\mu^{-1}+\mu\BI)}\big\|  \le \eop^{t@\nabs(\BA, \BE)}\]
for all $t\ge 0$, and all solutions of the DAE $\BE@\Bx'(t) = \BA\Bx(t)$ 
satisfy
\begin{equation} \label{eq:nabsub}
 {\|\Bx(t)\|\over \|\Bx(0)\|}  \le \eop^{t@\nabs(\BA,\BE)}.
\end{equation}
\end{theorem}

This bound suffers from a major limitation:  if $(\BA,\BE)$ 
is stable but $\nabs(\BA,\BE)>0$ (\new{as is true for the examples in \cref{fig:psa1}}),
\cref{eq:nabsub} fails to capture $\|\Bx(t)\| \to 0$ as $t\to\infty$.
To describe that convergence, suppose one can diagonalize
$\BG_\mu^{-1} = \BV \BLambda_\mu\BV^{-1}$,
so that
\begin{equation} \label{eq:evbnd}
 {\|\Bx(t)\| \over \|\Bx(0)\|} \le \|\eop^{t(\BG_\mu^{-1}+\mu\BI)} \| 
                                 \le \|\BV\|@\|\BV^{-1}\|\,\eop^{t@\alpha(\BA,\BE)},
\end{equation}
where $\alpha(\BA,\BE)$ is the spectral abscissa of $(\BA,\BE)$, 
i.e., the real part of the rightmost (finite) eigenvalue of $(\BA,\BE)$.\ \ 
If $(\BA,\BE)$ is stable, 
then $\alpha(\BA,\BE) <0$ and~\cref{eq:evbnd} describes $\|\Bx(t)\| \to 0$. 
However, $\|\BV\|@\|\BV^{-1}\|$ can be very large 
(or $\BG_\mu^{-1}$ may not be diagonalizable),
and this quantity is difficult to estimate when the pencil
has large dimension.
Pseudospectra give more flexible bounds that are 
\new{well-suited to approximation in the large-scale case
(as demonstrated in the next section)}.

\begin{theorem}
Let $(\BA,\BE)$ be a regular pencil with $\BA-\mu@\BE$ invertible.
For all  $\eps>0$ and $t>0$, 
\begin{equation} \label{eq:psabnd}
 \|\eop^{t(\BG_\mu^{-1}+\mu\BI)}\| 
     \ \le\ { L_\eps \eop^{t \alpha_\eps(\BA,\BE)} \over 2@\pi \eps}, 
\end{equation}
where $L_\eps$ is the contour length of a Jordan curve that contains 
$\PSA(\BA,\BB)$ in its interior.

For all $t\ge 0$, 
\begin{equation} \label{eq:kmt}
 \|\eop^{t(\BG_\mu^{-1}+\mu\BI)}\| \le \eop@@(n-d)@@{\cal K}(\BA,\BE),
\end{equation}
where $n-d$ is the dimension of $\BG_\mu$ and ${\cal K}(\BA,\BE)$ 
denotes the Kreiss constant.
\end{theorem}

Varying $\eps>0$ in~\cref{eq:psabnd} leads to a family of upper bounds:
as $\eps\downarrow 0$,
$\psabs(\BA,\BE)$ decreases monotonically to $\alpha(\BA,\BE)$
while $L_\eps/(2@\pi\eps)$ generally increases.  
The bound~\cref{eq:psabnd} is derived by crudely estimating the 
norm of the Dunford--Taylor integral~\cite[p.~44]{Kat76}
\begin{equation} \label{eq:dunford}
 \eop^{t(\BG_\mu^{-1}+\mu\BI)} 
  = {1\over 2\pi \new{@\iop}} 
     \int_{\Gamma_\eps} \eop^{t@z} (z\BI - (\BG_\mu^{-1}+\mu\BI))^{-1}\,\dop z,
\end{equation}
where $\Gamma_\eps$ is a finite union of Jordan curves enclosing $\PSA(\BA,\BE)$ 
in their collective interior.
When $L_\eps$ is large because $\Gamma_\eps$ must capture portions 
of $\PSA(\BA,\BE)$ far in the left-half plane, more careful estimates of the
integral~\cref{eq:dunford} could yield tighter bounds.

For stable $(\BA,\BE)$, since $\alpha(\BA,\BE) < 0$ one can
take $\eps>0$ sufficiently small that $\psabs(\BA,\BE)<0$.
For such $\eps$, \cref{eq:psabnd} implies $\|\Bx(t)\| \to 0$ as $t\to\infty$.
The leading constant $L_\eps/(2\pi \eps)$ then limits the extent
of transient growth.  
The bound~\cref{eq:kmt}, \new{known as} the Kreiss Matrix Theorem,
has a nontrivial proof with an interesting history behind the
dimension-dependent factor; see~\cite[chap.~18]{TE05}, \cite{WT94}.

\new{In summary}, 
any bound on $\|\eop^{t\BA}\|$ leads to a similar bound for DAEs
by simply replacing $\BA$ with $\BG_\mu^{-1}+\mu\BI$.
The sampling of bounds above is not meant to be exhaustive.
For example, one can obtain more refined (but complicated) bounds using 
pseudospectra~\cite[chap.~15]{TE05},
or by decomposing $\BG_\mu^{-1}$ using spectral projectors.
A rather different class of bounds involves the solution of 
an associated Lyapunov equation;
see, e.g., \cite[sect.~11.4]{God98}, \cite{Ves03a}, \cite[thm.~13.6]{Ves11}.

\section{Approximation of pseudospectra for large scale problems}
\label{sec:largescale}

For large $\BA$ and $\BE$, as often arise in
linear stability analysis problems derived 
from partial differential equations,
it is impractical to compute the sets $\PSA(\BA,\BE)$ in \cref{def:psa}.
For example, fluid dynamics applications give DAEs with coefficients of the form
\begin{equation} \label{eq:kkt}
 \BA = \left[\begin{array}{cc}
                \BK & \BB^*\! \\ \BB & \Bzero\!  
               \end{array}\right], \qquad
   \BE = \left[\begin{array}{cc}
                \BM & \Bzero \\ \Bzero & \Bzero  
               \end{array}\right],
\end{equation}
with $\BK\in \R^{n_v \times n_v}$ invertible, 
$\BB \in \R^{n_p\times n_v}$ full rank, 
and $\BM \in \R^{n_v \times n_v}$ Hermitian positive definite,
\new{for $n_v \ge n_p$}.
The pencil $(\BA,\BE)$ has $n_v-n_p$ finite eigenvalues
and $2n_p$ infinite eigenvalues (associated with $n_p$ Jordan
blocks of size $2\times 2$), so the corresponding DAE 
has index~2; see~\cite{CGS94} for a discussion of this
eigenvalue problem.
Given this spectral structure, 
$\BG_\mu \in \C^{(n_v-n_p)\times(n_v-n_p)}$;
in engineering computations $n_v-n_p$ can easily be $10^4$ or much larger.
Our proposed definition of pseudospectra will only be useful
if there is a practical way to compute approximations
that require little effort beyond the standard eigenvalue computation 
already used for linear stability analysis.

\cool{Here we continue using the 2-norm; the technique is extended
to alternative norms in~\cref{sec:altnorms}.}
Wright and Trefethen proposed a technique for approximating 
conventional pseudospectra by restricting the matrix to an
invariant subspace~\cite{WT01b} computed using ARPACK~\cite{LSY98}
(perhaps via MATLAB's {\tt eigs} interface).  
This approach provides interior 
estimates of the pseudospectra; i.e., if the columns of $\BV\in\C^{n\times k}$ 
form an orthonormal basis for a $k$-dimensional invariant subspace of $\BA$,
then for all $\eps>0$,
\begin{equation} \label{eq:psabound}
 \PSA(\BV^*\!\BA\BV)\subseteq \PSA(\BA).
\end{equation}
(The EigTool software offers a modified projection method, 
where the invariant subspace
is augmented by a Krylov subspace, \new{with interior bounds obtained
from pseudospectra of \emph{rectangular} Hessenberg matrices}~\cite{Wri02a,WT01b}.)
If the invariant subspace corresponds to
all eigenvalues in some region of the complex plane
(e.g., the rightmost eigenvalues),
then the approximation~\cref{eq:psabound}
is typically quite accurate near those eigenvalues.
(See~\cite[chap.~40]{TE05}, which also explains
when this approximation fails to be accurate.) 
The matrix $\BV^*\!\BA\BV \in \C^{k\times k}$ is generally
much smaller than $\BA$, so its pseudospectra can be 
computed using standard dense techniques~\cite{Tre99a}
in a fraction of the time it took to compute $\BV$.
Thus approximate pseudospectra can be generated 
as a simple byproduct of a large-scale eigenvalue
computation, providing a simple way to perform a 
pseudospectral sensitivity analysis.

We seek a similar approximation strategy for the
pseudospectra of the matrix pencil, $\PSA(\BA,\BE)$.
To assess the asymptotic stability of solutions of 
the DAE~\cref{eq:dae},
one seeks the rightmost (finite) eigenvalues of the 
pencil $(\BA,\BE)$; these are typically found by computing the
largest-magnitude eigenvalues of the shift-invert 
transformation $(\BA-\mu@\BE)^{-1}\BE$ or Cayley 
transformation $(\BA-\mu_1\BE)^{-1}(\BA-\mu_2@\BE)$;
see, e.g., \cite{MSR94}.

Suppose that for $\mu\in\C$, the matrix $\BA-\mu\BE$ is invertible,
and let the columns of $\BV\in\C^{n\times k}$ give an 
orthonormal basis for a $k$-dimensional invariant
subspace of $(\BA,\BE)$ associated with finite eigenvalues.
(Equivalently, $\Ran(\BV)$ is an invariant subspace
of $\BE_\nu := (\BA-\nu\BE)^{-1}\BE$ associated with
nonzero eigenvalues for any $\nu\in\C$ for which 
$\BA-\nu\BE$ is invertible, following essentially the
same argument that showed $\mu$-independence of
\cref{def:psa}.  
Thus $\BV$ can be computed using any desired shift-invert 
transformation.)
Now $\SPEC(\BV^*\Emu\BV) \subseteq \SPEC(\Emu)$; in particular,
consider the Schur factorization of the $k\times k$ matrix
\[ \BV^*\Emu\!\BV = \BU@\Ghat\BU^*,\]
where $\Ghat\in\C^{k\times k}$ is an invertible upper-triangular matrix
with $\SPEC(\Ghat)\subseteq \SPEC(\Emu)$ and 
$\BU\in\Ckk$ is unitary.
This decomposition is a partial Schur factorization of
$\Emu$: since the eigenvalues can be ordered arbitrarily on the 
diagonal of the Schur factor, we can compute some unitary 
$[\BQ\ \BQ_\perp] \in \Cnn$
such that
\[ 
 \Emu = \left[\begin{array}{cc} \BQ & \BQ_\perp \end{array}\right]
          \left[\begin{array}{cc} \BG & \BD \\ \Bzero & \BN \end{array}\right]
          \left[\begin{array}{c} \BQ^* \\ \BQ_\perp^* \end{array}\right]
\]
with 
\begin{equation} \label{eq:G}
 \BG = \left[\begin{array}{cc}  \Ghat & \BX \\ \Bzero & \wt{\BG} \end{array}\right]
\end{equation}
for $\Ghat \in \C^{k\times k}$ and $\wt{\BG}\in\C^{(n-d-k)\times (n-d-k)}$
both invertible.
To compute $\PSA(\BA,\BE)$ in \cref{def:psa}, we must
compute level sets of $\|((z-\mu)\BI -\BG^{-1})^{-1}\|$.
Note that 
\[ \BG^{-1} = \left[\begin{array}{cc}
             \Ghat^{-1} & -\Ghat^{-1}\BX\wt{\BG}^{-1} \\[.25em]
             \Bzero     & \wt{\BG}^{-1}
              \end{array}\right],\]
and 
\[ ((z-\mu)\BI -\BG^{-1})^{-1}
           = \left[\begin{array}{cc}
             ((z-\mu)\BI - \Ghat^{-1})^{-1} & \BXi \\[.25em]
             \Bzero     & ((z-\mu)\BI - \wt{\BG}^{-1})^{-1}
              \end{array}\right]\]
for $\BXi := -((z-\mu)\BI-\Ghat^{-1})^{-1}\Ghat^{-1}\BX\wt{\BG}^{-1}((z-\mu)\BI-\wt{\BG}^{-1})^{-1}$.
The \cool{2-}norm of the $(1,1)$ block of \cool{$((z-\mu)\BI-\BG^{-1})^{-1}$}
cannot exceed the \cool{2-}norm of the entire matrix, so 
\begin{align*} \label{eq:resbnd}
  \| ((z-\mu)\BI -\BG^{-1})^{-1}\| 
    &\ge \| ((z-\mu)\BI -\Ghat^{-1})^{-1}\| \\[.25em]
    &=   \| (z\BI - (\Ghat^{-1}+\mu\BI))^{-1}\|. \nonumber 
\end{align*}  
Applying this bound to \cref{def:psa} shows that the
computed invariant subspace gives an \emph{interior bound} 
on the pseudospectra of $(\BA,\BE)$.
(For simplicity of formulation, we omit the unitary similarity
transformation with $\BU$ from the definition of $\Ghat$, 
as it does not alter the pseudospectra.)

\begin{theorem} \label{thm:psabnd}
Let the columns of $\BV\in\C^{n\times k}$ form an orthonormal basis
for a $k$-dimensional invariant subspace of $(\BA,\BE)$ associated 
with finite eigenvalues, and let $\mu\in\C$ be any number for which
$\BA-\mu\BE$ is invertible.
Then for all $\eps>0$,
\[ \PSA(\Ghat^{-1}+\mu\BI) \subseteq \PSA(\BA,\BE),\]
where $\Ghat = \BV^*(\BA-\mu\BE)^{-1}\BE\BV$.
\end{theorem}

\medskip
This theorem implies that lower bounds on $\PSA(\BA,\BE)$ can be obtained as
a byproduct of the usual eigenvalue calculation performed for
linear stability analysis.  Two caveats are in order.  
(1)~To obtain pseudospectral estimates in the norm most relevant for 
the physical problem, one should first transform $\BA$ and $\BE$ as described
in \cref{sec:norms}, so that the \new{2}-norm on $\Cn$ gives
an accurate measure of the physically-motivated norm.  The basis vectors
for the invariant subspace in $\BV$ are thus orthogonal in the Euclidean
\new{inner product}.
(2)~To accurately approximate $\PSA(\BA,\BE)$, one often needs a large
invariant subspace, i.e., $k$ might be taken larger than one would 
use if only computing the rightmost eigenvalue.
However, larger subspaces bolster one's confidence that a rightmost 
eigenvalue with large imaginary part has not been missed, and further
reveal the role of subordinate eigenvalues on the transient behavior.
The next section shows several illustrations for problems from fluid
dynamics.

The accuracy of the approximation in \cref{thm:psabnd} depends on
several factors, such as the location of the computed eigenvalues, 
the dimension of the associated invariant subspace, 
and the angle between that subspace and the complementary 
invariant subspace associated with the other \emph{finite} eigenvalues
(related to the matrix $\BX$ in~\cref{eq:G}).
We cannot expect the approximation to be accurate throughout $\C$,
particularly when $k\ll n-d$.
Rather, we hope it is accurate in a region of $\C$ most relevant
to the application at hand.  
For example, for linear stability analysis of a continuous-time 
dynamical system, we hope the approximation
$\PSA(\Ghat^{-1}+\mu\BI) \approx \PSA(\BA,\BE)$ is accurate 
in the intersection of the right-half plane with $W(\BG^{-1}+\mu\BI)$,
which will lead to accurate estimates of the positive values of $\alpha_\eps(\BA,\BE)$.
In any case, since \cref{thm:psabnd} gives interior estimates, we always
have $\alpha_\eps(\Ghat^{-1}+\mu\BI) \le \alpha_\eps(\BA,\BE)$. 
Thus \cref{thm:explow} implies the following lower bound on transient growth.

\begin{corollary} \label{cor:approx_etA}
Using the notation of \cref{thm:psabnd},
for any $\eps>0$
there exists some initial condition $\Bx(0)\in\Ran(\BQ_\mu)$ 
such that the solution $\Bx(t)$ to~$\cref{eq:dae}$ satisfies
\[ \sup_{t\ge 0} {\|\Bx(t)\| \over \|\Bx(0)\|} \ge {\alpha_\eps(\Ghat^{-1})+\mu \over \eps}.\]
\end{corollary}

\subsection{Alternative Norms} \label{sec:altnorms}

Suppose the matrix pencil is derived from a physical problem that is 
associated with some domain-specific inner product. 
Practical eigenvalue computations for linear stability analysis
usually make no special effort to compute with this physically 
relevant inner product:
the inner product does not affect the eigenvalues of the matrix pencil,
and use of a different inner product would incur additional arithmetic
beyond that needed for the standard 2-norm calculation.

Suppose we have a matrix $\BV\in\C^{n\times k}$ whose columns form
a basis for an invariant subspace associated with nonzero eigenvalues
of $(\BA-\mu@\BE)^{-1}\BE$ that is orthonormal in the 2-norm, 
so there exists some $\wh{\BG}\in\C^{k\times k}$ such that
\begin{equation} \label{eq:invsubsp}
 (\BA-\mu@\BE)^{-1}\BE\BV = \BV\wh{\BG}.
\end{equation}
Using the notation of \cref{sec:norms}, we wish to 
approximate $\sigma_{\eps,\BH}(\BA,\BE)$, where $\BH$ is a
positive definite matrix with the factorization $\BH = \BR^*\BR$.
To approximate these $\BH$-norm pseudospectra using the approach 
outlined earlier in this section, it will suffice to transform $\BV$ to
obtain a 2-norm orthonormal basis for the corresponding invariant
subspace of $\BR(\BA-\mu@\BE)^{-1}\BE\BR^{-1}$ (see~\cref{eq:EmuHdef}).

Now~\cref{eq:invsubsp} is equivalent to
\[ \BR(\BA-\mu@\BE)^{-1}\BE\BR^{-1}\BR\BV = \BR\BV\wh{\BG}.\]
Compute an economy-sized QR factorization $\BR\BV = \BZ\BS$, 
so $\BZ^*\BZ=\BI\in\C^{k\times k}$ and 
\[ \BZ^*\big(\BR(\BA-\mu@\BE)^{-1}\BE\BR^{-1}\big)\BZ = \BS\wh{\BG}\BS^{-1}.\]
Using the same arguments behind \cref{thm:psabnd} and \cref{cor:approx_etA},
we have 
\begin{equation} \label{eq:approxH}
 \sigma_{\eps,2}(\BS\wh{\BG}^{-1}\BS^{-1}+\mu@\BI) \subseteq \sigma_{\eps,\BH}(\BA,\BE) 
\end{equation}
and
\begin{equation} \label{eq:etAH}
 \sup_{t\ge 0} {\|\Bx(t)\|_\BH \over \|\Bx(0)\|_\BH} 
      \ge {\alpha_{\eps,2}(\BS\Ghat^{-1}\BS^{-1})+\mu \over \eps}.
\end{equation}
Thus, pseudospectra can be readily approximated in physically relevant norms
using the invariant subspace $\BV$ deriving from a standard 2-norm linear
stability analysis.  

Related ideas for approximating standard pseudospectra in weighted
norms are described by Astudillo and Castillo~\cite{AC13}.
We also note that the new reduced basis techniques for standard pseudospectra 
of Sirkovi\'c~\cite{Sir} also hold great promise for estimating $\PSA(\BA,\BE)$.

At the end of \cref{sec:Hnorm2}, we summarized how one can compute
$\sigma_{\eps,\BH}(\BA,\BE)$ for small- or medium-scale problems.
Here we provide a similar summary for approximating
$\sigma_{\eps,\BH}(\BA,\BE)$ in the large-scale case, given $\BH=\BR^*\BR$.
\begin{enumerate}
\item Compute the $k$-dimensional invariant subspace of $(\BA-\mu\BE)^{-1}\BE$
associated with the eigenvalues of most relevance to the application.
Let the columns of $\BV\in\C^{n\times k}$ give a basis for this subspace 
that is orthonormal in the 2-norm, and let $\wh{\BG}\in\Ckk$ be the
generalized Rayleigh quotient given in~\cref{eq:invsubsp}.
\item Compute the economy-sized QR factorization $\BR\BV = \BZ\BS$.
\item Compute the lower bound $\sigma_{\eps,2}(\BS\wh{\BG}^{-1}\BS^{-1} + \mu@\BI)\subseteq \sigma_{\eps,\BH}(\BA,\BE)$.
\end{enumerate}

\section{Computational examples}
\label{sec:examples}

\Cref{fig:psa1} showed pseudospectra for two matrix pencils of size $n=3$. 
In this section we study pseudospectra for much larger problems 
that arise from linear stability analysis for several incompressible fluid flows 
in two physical dimensions.  These examples were generated using the IFISS 
software package~\cite{ERS14}.%
\footnote{We are grateful to Howard Elman for considerable guidance with this
software, and for sharing code to generate and extract the requisite matrices 
from within IFISS.} 
(Note the recent work of Emmrich and Mehrmann~\cite{EM13}, which compares
the spatial discretization approach used here to direct analysis of 
infinite dimensional fluid DAEs.)
Given a domain $\Omega \subset \R^2$, 
the velocity field $\Su: \Omega \times (0,\infty) \to \R^2$ 
and pressure field $\Sp:\Omega \times (0,\infty) \to \R$ 
satisfy the incompressible Navier--Stokes equations 
\begin{align*} 
   \Su_t(\Sx,t) &= -\nu \Delta \Su(\Sx,t) 
                    + \Su(\Sx,t) \cdot \nabla \Su(\Sx,t) + \nabla \Sp(\Sx,t) \\
             0  &= \nabla \cdot \Su(\Sx,t)
\end{align*}
for $\Sx\in \Omega \subset \R^2$, with proper boundary conditions for the flow.  
Here $\nu>0$ denotes the \emph{viscosity}, which is inversely proportional 
to the Reynolds number.
We first seek a steady-state solution $(\widehat{\Su}(\Sx), \widehat{\Sp}(\Sx))$ for which
$ -\nu \Delta \wh{\Su}(\Sx) + \wh{\Su}(\Sx) \cdot \nabla\wh{\Su}(\Sx) + \nabla \wh{\Sp}(\Sx)=\Bzero$
and $\nabla \cdot \wh{\Su}(\Sx) = 0$.
Is this stationary solution stable when subjected to small perturbations?
Linear stability analysis (see, e.g., \cite[chap.~15--16]{Gun89}) inserts
 $\Su(\Sx,t) := \wh{\Su}(\Sx) + \Sw(\Sx,t)$ and $\Sp(\Sx,t) := \wh{\Sp}(\Sx) + \Ss(\Sx,t)$
\new{(with small $\|\Sw\|$ and $\|\Ss\|$)}
into the incompressible Navier--Stokes equations and neglects the 
quadratic term $\Sw \cdot \nabla \Sw$ (since $\|\Sw\|\ll 1$) to 
approximate evolution of the perturbation as
\begin{align*} 
   \Sw_t(\Sx,t) &= -\nu \Delta \Sw(\Sx,t) 
                    + \wh{\Su}(\Sx) \cdot \nabla \Sw(\Sx,t) + \Sw(\Sx,t)\cdot \nabla \wh{\Su}(\Sx) + \nabla \Ss(\Sx,t) \\
             0  &= \nabla \cdot \Sw(\Sx,t).
\end{align*}
Common finite element discretizations of this equation yield a DAE of the form
\begin{equation} \label{eq:nsdae}
 \left[\begin{array}{cc} \BM & \Bzero \\ \Bzero & \Bzero \end{array}\right]
   \left[\begin{array}{cc} \!\!\Bw'(t)\!\! \\ \!\!\Bs'(t)\!\! \end{array}\right]
   = \left[\begin{array}{cc} \BK & \BB^* \\ \BB & \Bzero \end{array}\right] 
   \left[\begin{array}{cc} \!\!\Bw(t)\!\! \\ \!\!\Bs(t)\!\! \end{array}\right],
\end{equation}
where $\BM, \BK \in \R^{n_v\times n_v}$ are invertible and 
$\BB\in \R^{n_p \times n_v}$ has full rank.
(Here $n_v$ and $n_p$ denote the number of discretized velocity and pressure
variables, with $n_v > 2@n_p$: $\Bw(t) \in \R^{n_v}$, $\Bs(t) \in \R^{n_p}$.)
Spectral properties of the associated pencil $(\BA,\BE)$ are discussed in~\cite{CGS94}.
The structure ensures that $(\BA,\BE)$ has an infinite eigenvalue of multiplicity 
$2@n_p$ associated with $n_p$ Jordan blocks, each of dimension~2.  
Hence, the DAE has index~2, and in the notation of~\cref{eq:schur},
we know \emph{a priori} that the block $\BN_\mu$ has dimension $d= 2@n_p$.

It is customary to measure the perturbations 
$\Sw$ and $\Ss$ via 
\begin{align*}
  |\Sw(\cdot,t)|_{H_1} 
     &= \bigg( \int_\Omega \|\nabla {\sf w}_1(\Sx,t)\|^2 
                                         + \|\nabla {\sf w}_2(\Sx,t)\|^2\,\dop \Sx\bigg)^{1/2}, \\[.25em]
   \|\Ss(\cdot,t)\|_{L_2} 
     &= \bigg( \int_\Omega |\Ss(\Sx,t)|^2 \,\dop \Sx\bigg)^{1/2},
\end{align*}
where the norms on the right-hand side of the definition of $|\Sw(\Bx,t)|_{H_1}$ 
are standard vector 2-norms in $\R^2$; 
see, e.g., \cite[sect.~8.4]{ESW14}, \cite[sect.~IV.2]{GR86}.
We thus analyze the discretization~\cref{eq:nsdae} using a discrete approximation
to the norm
\begin{equation} \label{eq:flunorm}
 \left\|\left[\begin{array}{cc} \Sw(\cdot,t) \\ \Ss(\cdot,t) \end{array}\right]\right\| := \bigg(|\Sw(\cdot,t)|_{H_1}^2 + \|\Ss(\cdot,t)\|_{L_2}^2\bigg)^{1/2}.
\end{equation}
All our examples use a uniform grid with $Q_2$--$Q_1$ finite elements~\cite{ERS14}.

\subsection{Backward facing step}
\begin{figure}
\begin{center}
\includegraphics[scale=0.7]{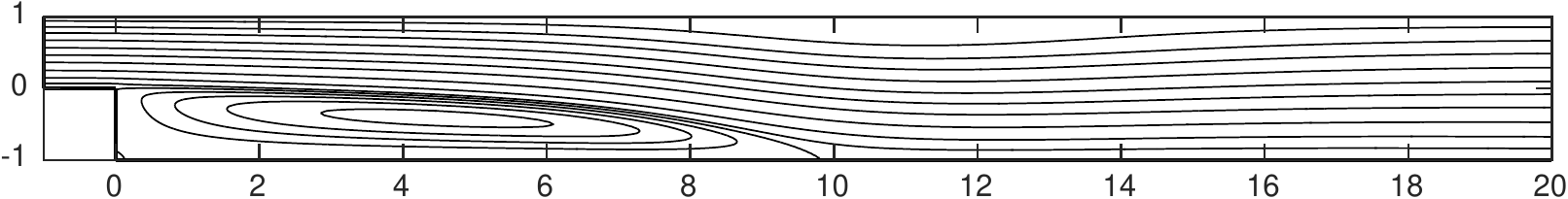}
\begin{picture}(0,0)
\put(-321.6,7){\color{white}{\rule{16pt}{15.65pt}}}  
\end{picture}
\end{center}

\vspace*{-7pt}
\caption{\label{fig:bfs}
Some (nonuniform) streamlines for the steady-state solution of the backward facing step problem
with viscosity $\nu=1/400$ with grid parameter ${\tt nc } = 6$ ($n_v = 85442$, $n_p = 10865$).
}
\end{figure}
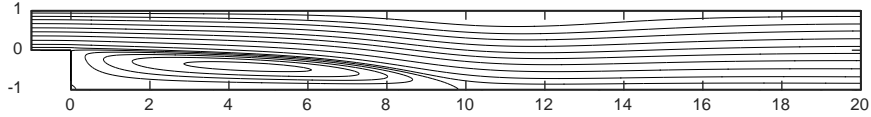

Our first example is the well-studied case of flow over a backward facing step; see, e.g., \cite{Gre93}.
Flow enters through the leftmost part of boundary and exits out the right end.
The step should be sufficiently long to resolve a dip in the streamlines 
near the top wall that moves further downstream as the viscosity $\nu$ decreases;
see \cref{fig:bfs}.
Indeed, to obtain satisfactory steady state flows as $\nu$ decreases, one must 
(a)~increase step length; (b)~refine the grid; (c)~exercise greater care with
the Picard and Newton nonlinear iterations used to find the steady state.
For all values of $\nu$ we have studied, the linearization is eigenvalue stable, 
i.e., all finite eigenvalues of $(\BA,\BE)$ are in the left-half plane, 
though \cool{the spectral abscissa approaches zero as $\nu$ decreases}.

After using IFISS to find the steady state flow for a given flow configuration, 
we approximate the pseudospectra of the pencil in~\cref{eq:nsdae} 
as described in \cref{thm:psabnd}:
use the {\tt eigs} command to compute the invariant subspace 
associated with the largest magnitude eigenvalues of $\BA^{-1}\BE$
(all calculations in this section use $\mu=0.25$),
orthonormalize these eigenvectors to obtain some $\BV\in\C^{n\times k}$,
and compute
$\PSA(\widehat{\BG}^{-1}+\mu\BI)\subseteq \PSA(\BA,\BE)$
for $\widehat{\BG} = \BV^*(\BAs-\mu\BE)^{-1}\BE\BV$
(in a discretization of the norm~\cref{eq:flunorm}.)

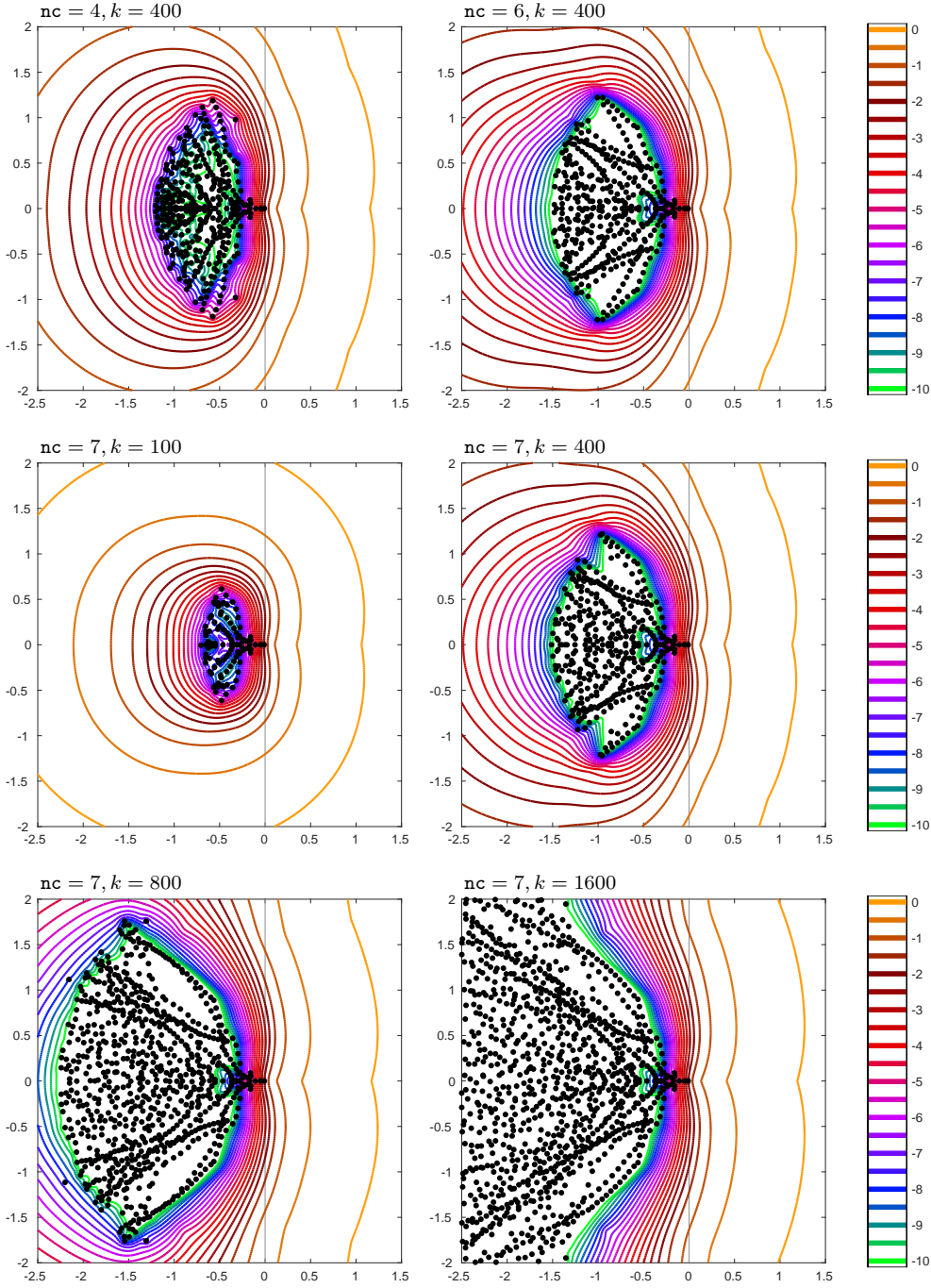
\begin{figure}[p]
\begin{picture}(0,500)(0,0)
\put(0,350){\includegraphics[scale=0.42]{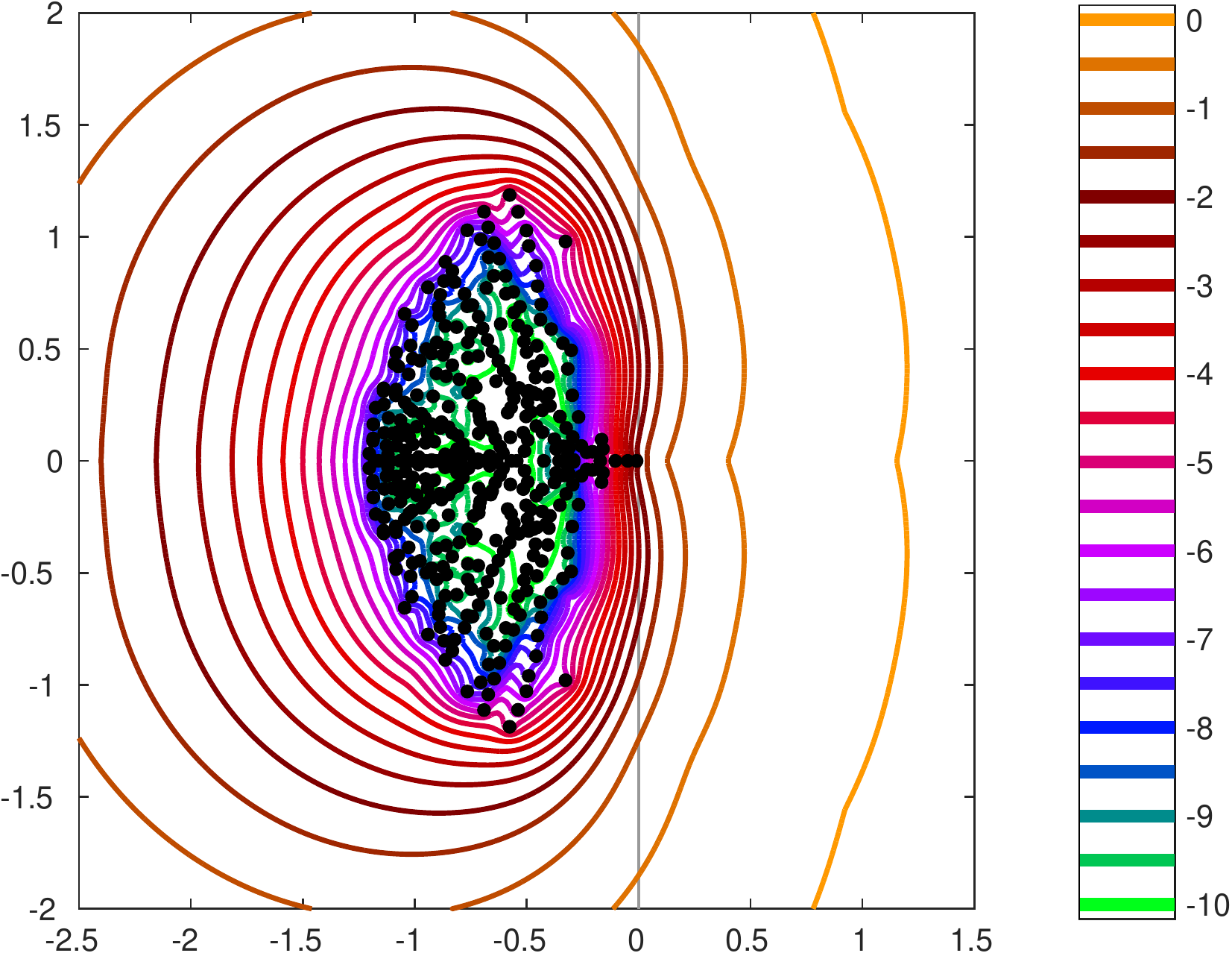}}
\put(170,350){\includegraphics[scale=0.42]{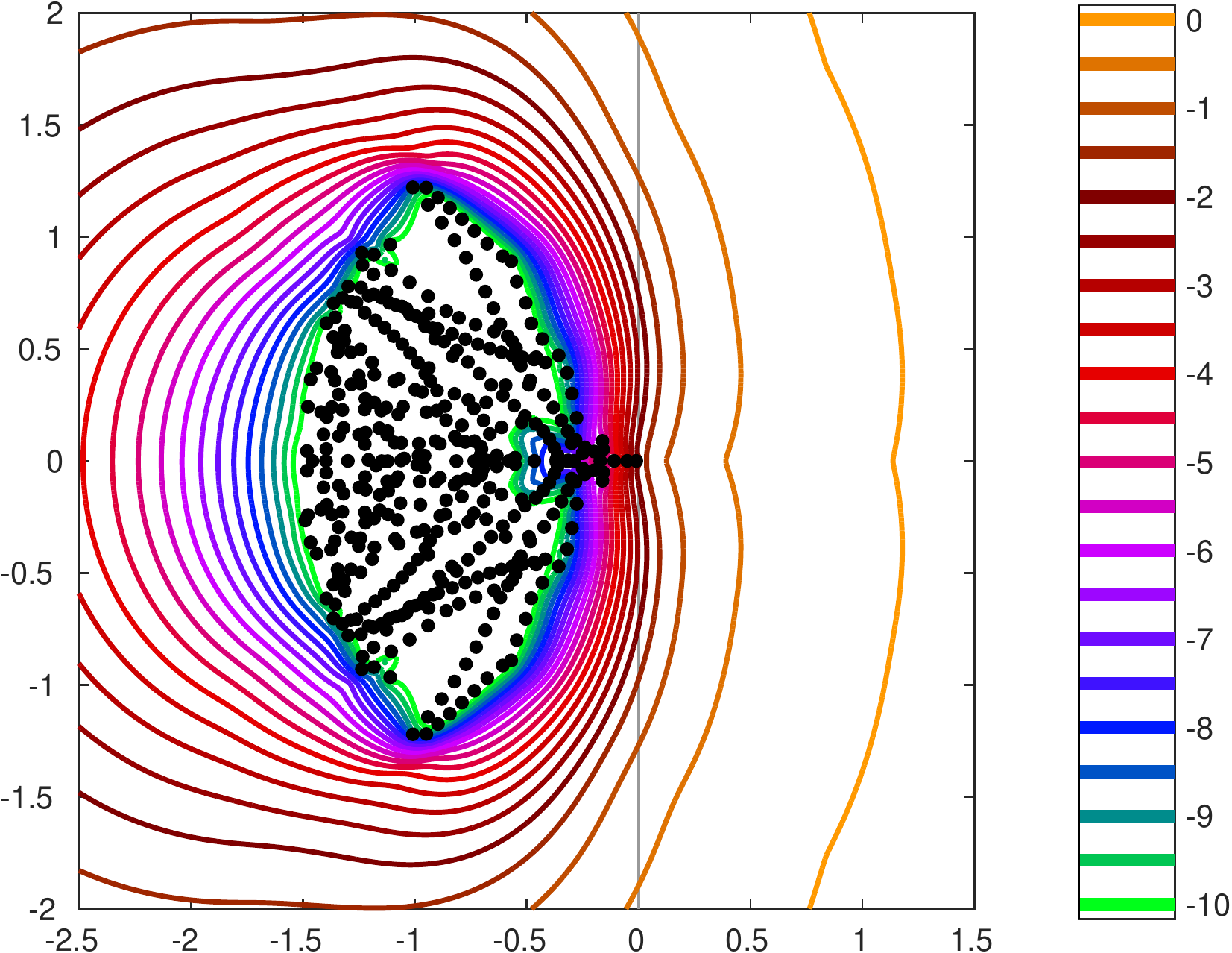}}
\put(0,175){\includegraphics[scale=0.42]{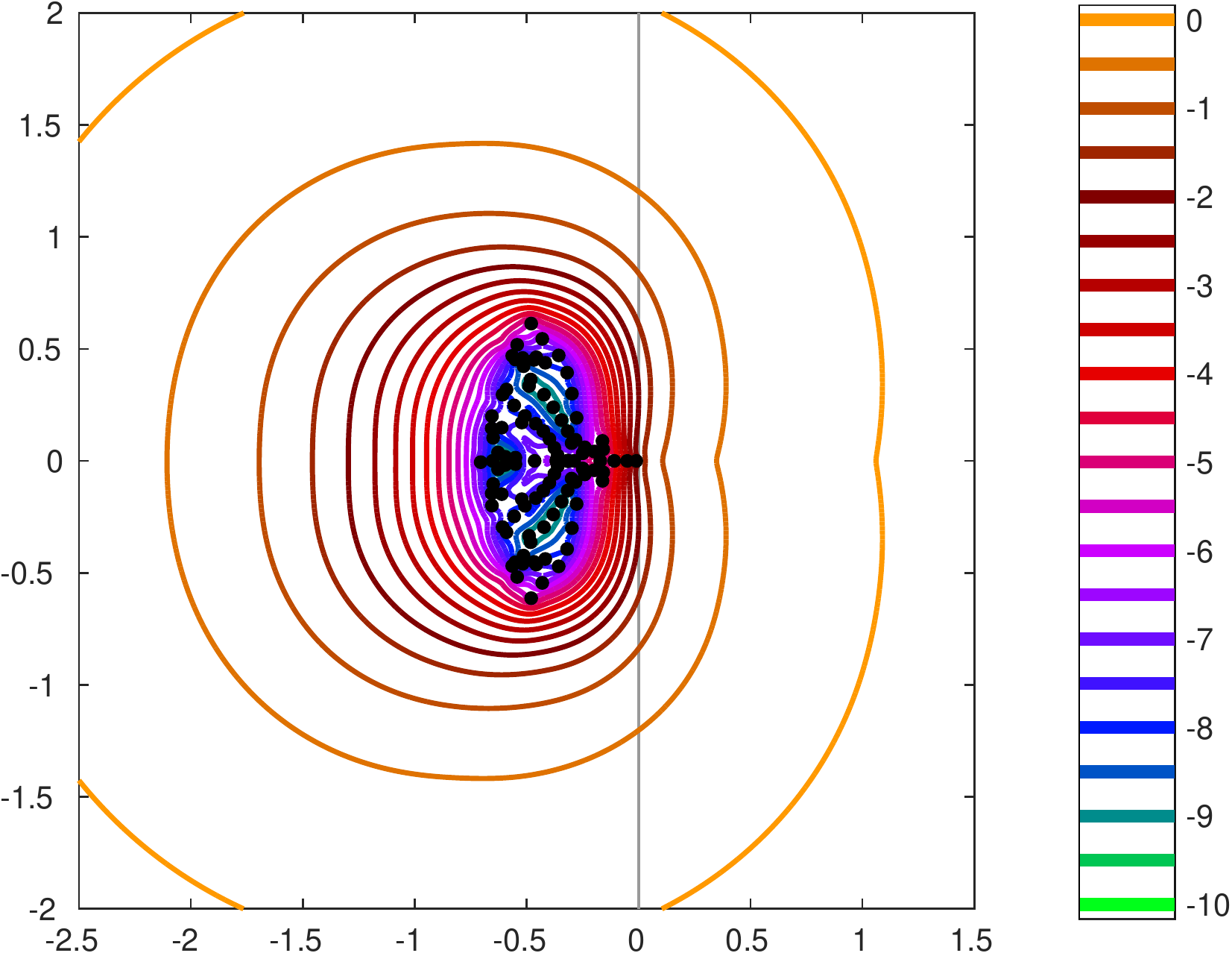}} 
\put(170,175){\includegraphics[scale=0.42]{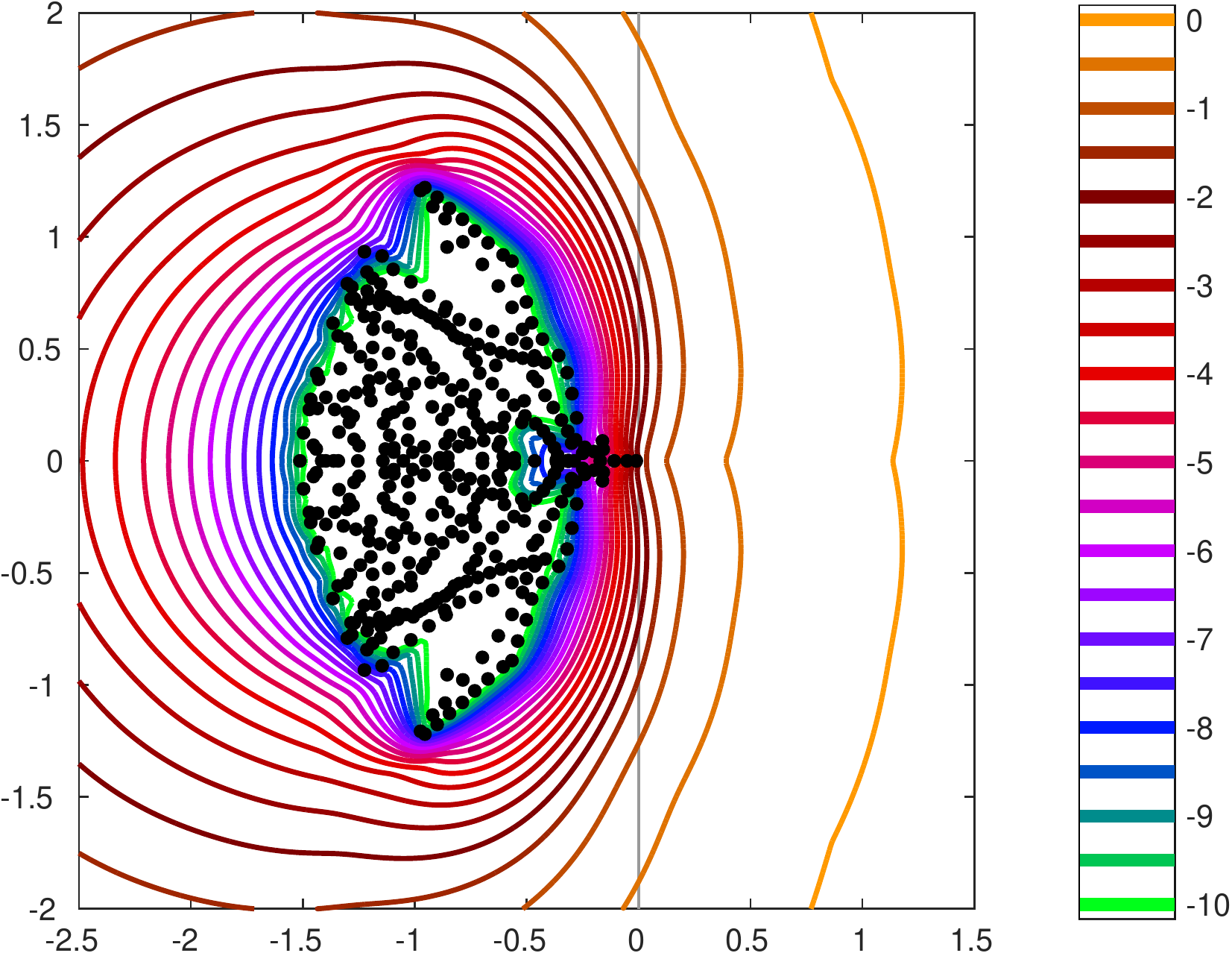}}
\put(0,0){\includegraphics[scale=0.42]{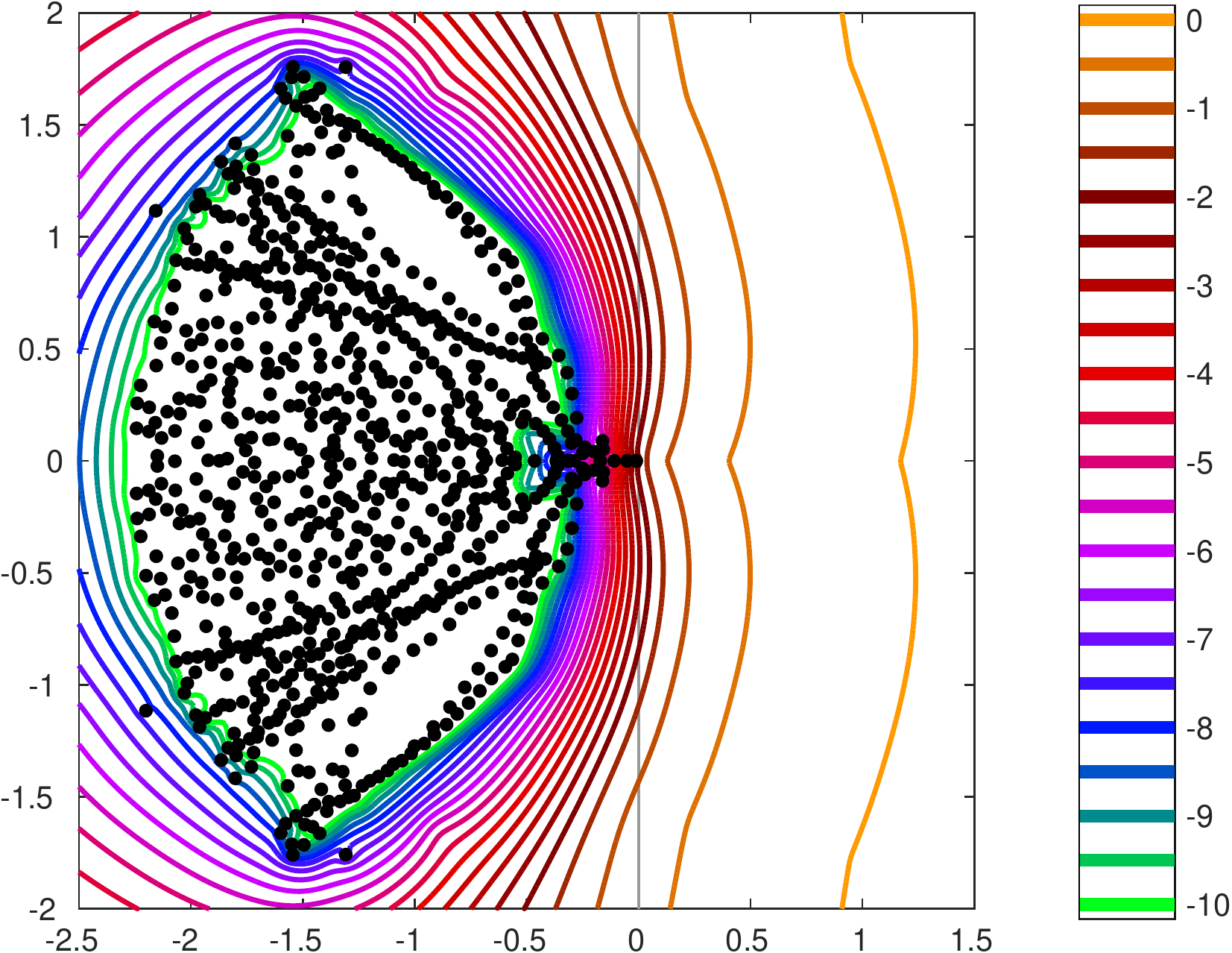}} 
\put(170,0){\includegraphics[scale=0.42]{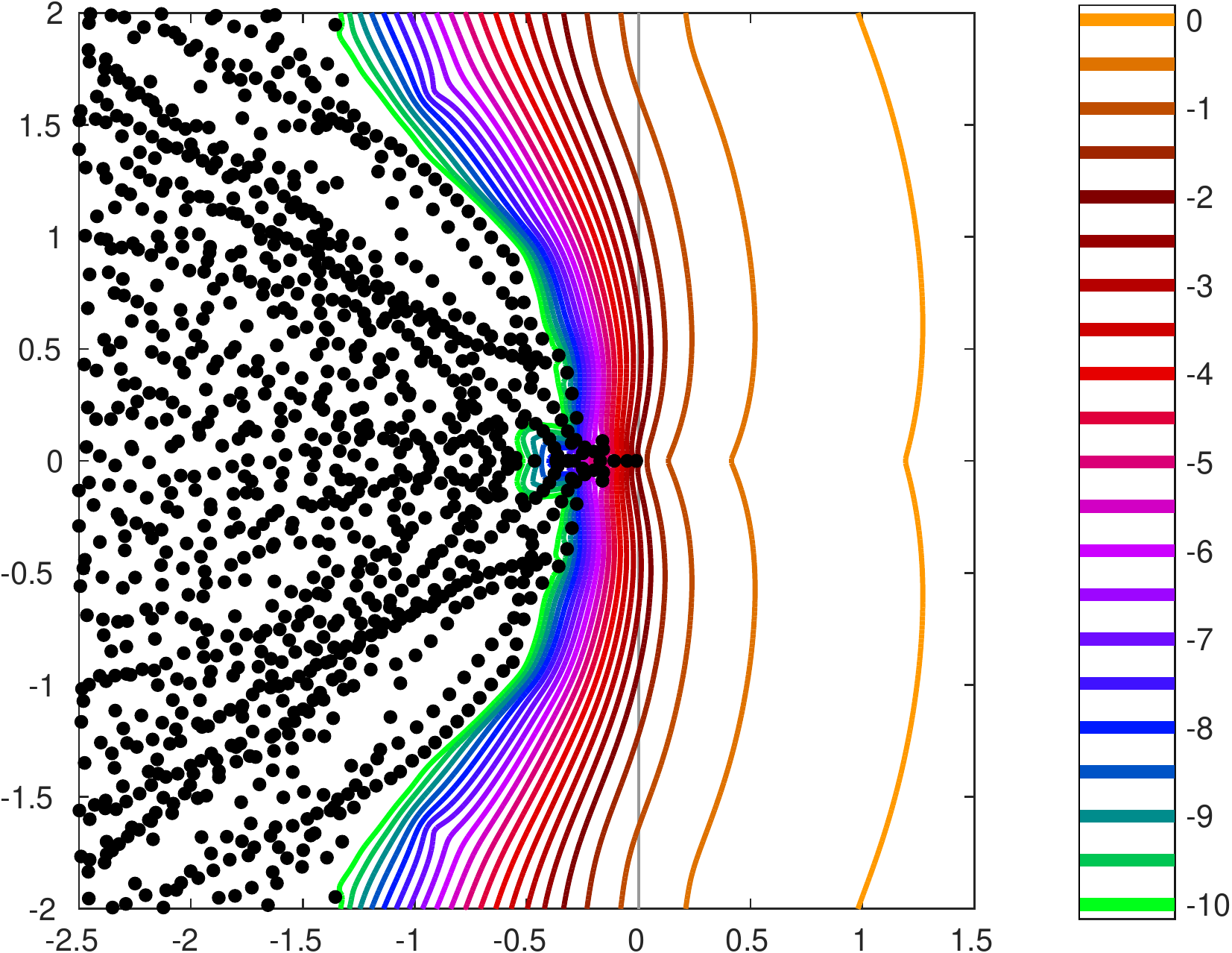}}
\put(14,507){\footnotesize ${\tt nc}=4, k=400$}
\put(184,507){\footnotesize ${\tt nc}=6, k=400$}
\put(14,332){\footnotesize ${\tt nc}=7, k=100$}
\put(184,332){\footnotesize ${\tt nc}=7, k=400$}
\put(14,157){\footnotesize ${\tt nc}=7, k=800$}
\put(184,157){\footnotesize ${\tt nc}=7, k=1600$}
\end{picture}

\caption{\label{fig:bfspsa}
Approximations of the pseudospectra $\PSA(\BA,\BE)$ 
for $\eps=10^{0}, 10^{-0.5}, \ldots, 10^{-10}$
for the backward facing step with viscosity $\nu=1/400$.
The top two plots use projection onto an invariant subspace of dimension $k=400$
for discretizations of dimension $n=6{,}367$ (${\tt nc}=4$) and 
       $n=96{,}307$ (${\tt nc}=6$).
The bottom four plots project a discretization of size
$n=381{,}539$ (${\tt nc}=7$) onto subspaces of dimension 
$k=100$, $400$, $800$, and $1600$.
The labels on the color bar show $\log_{10}\eps$, 
so, e.g., the orange contour on the right corresponds to
$\eps = 10^0$.
}
\end{figure}

\Cref{fig:bfspsa} shows various approximations to $\PSA(\BA,\BE)$ 
to illustrate several issues that arise when computing pseudospectra 
of large problems.
Three of the plots show estimates to $\PSA(\BA,\BE)$ using projection 
onto computed invariant subspaces of dimension $k=400$, based on
original \new{IFISS} discretizations of size 
${\tt nc}=4$ ($n=6{,}367$), 
${\tt nc}=6$ ($n=96{,}307$), and 
${\tt nc}=7$ ($n=381{,}539$).
The results change noticeably from ${\tt nc}=4$ to ${\tt nc}=6$,
but much less so from ${\tt nc}=6$ to ${\tt nc}=7$.

Four of the plots fix ${\tt nc}=7$, but project onto computed
invariant subspaces of varying dimension: $k=100$, $400$,
$800$, and $1600$.  To gain insight into the physical problem,
one cares about the extent of the pseudospectra into the 
right-half plane.  For example, since the boundary of the 
$\eps=10^0$ pseudospectrum extends beyond~$10^0$ in the real direction,
\cref{cor:approx_etA} ensures that, for some valid initial 
conditions, the differential algebraic equation will experience
transient growth.  
Note that even though the rightmost eigenvalue is real, 
the rightmost extent of the pseudospectra in these plots 
occurs at non-real values.

\begin{figure}[p]
\begin{picture}(0,500)(0,0) 
\put(0,350){\includegraphics[scale=0.42]{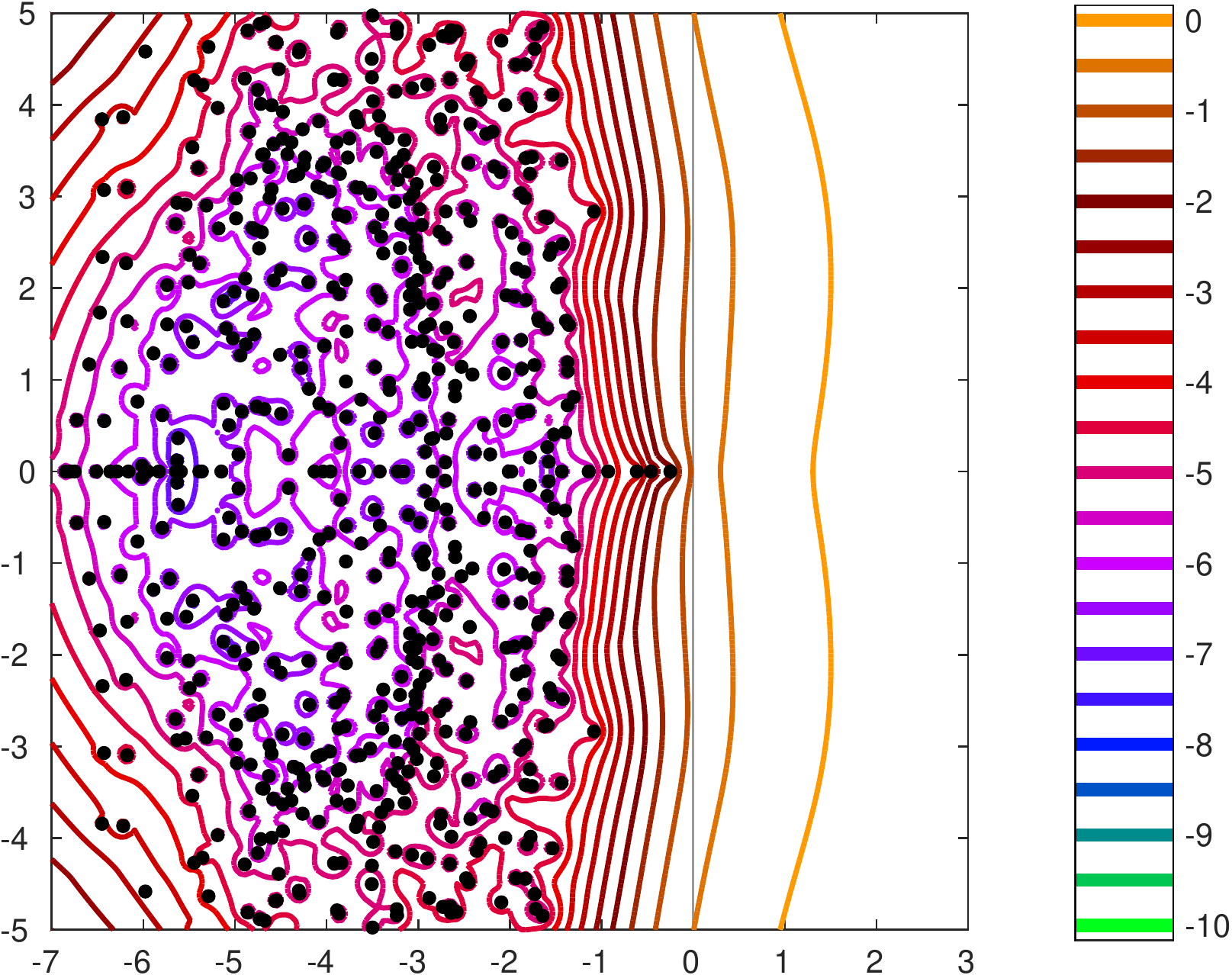}} 
\put(170,350){\includegraphics[scale=0.42]{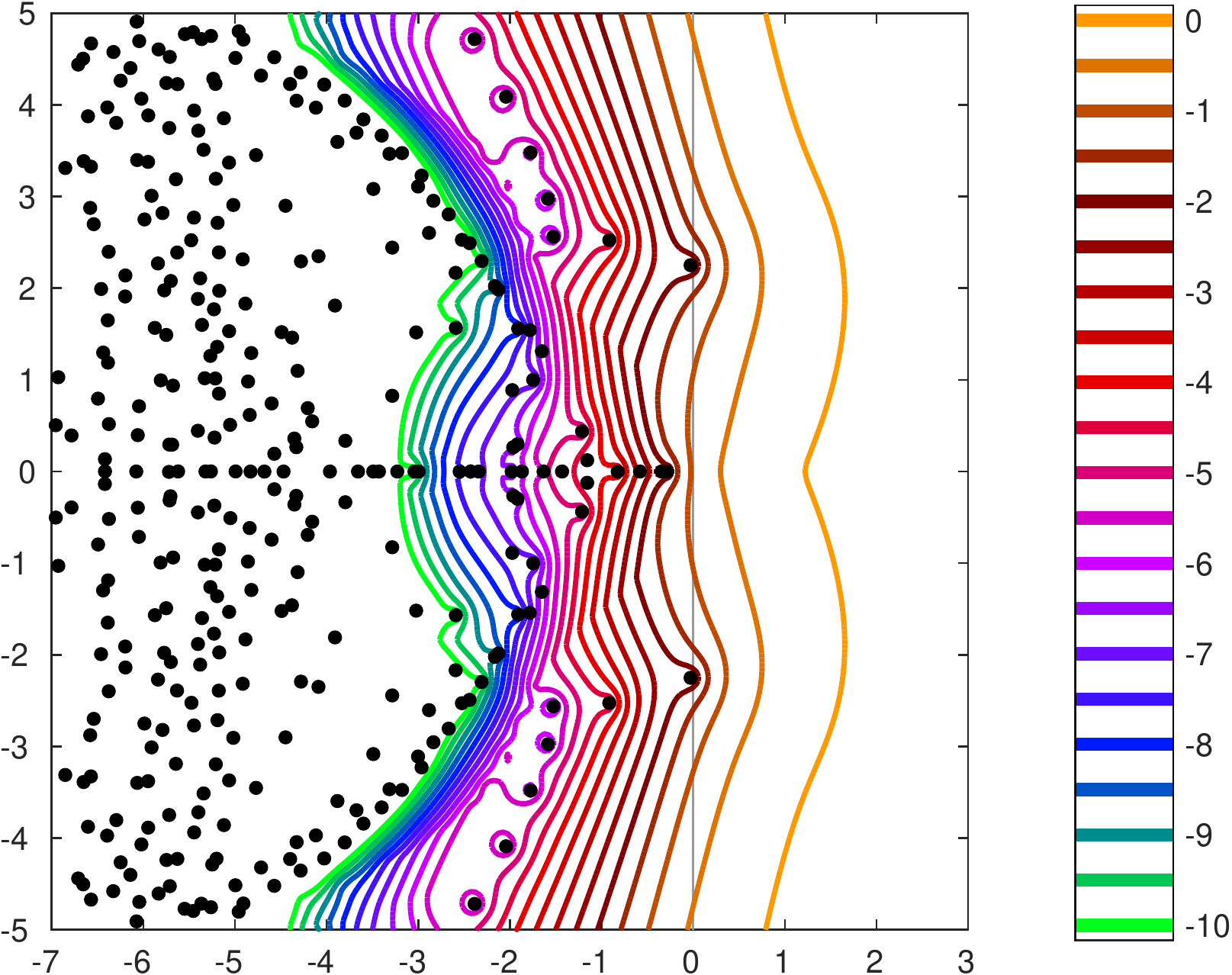}}
\put(0,175){\includegraphics[scale=0.42]{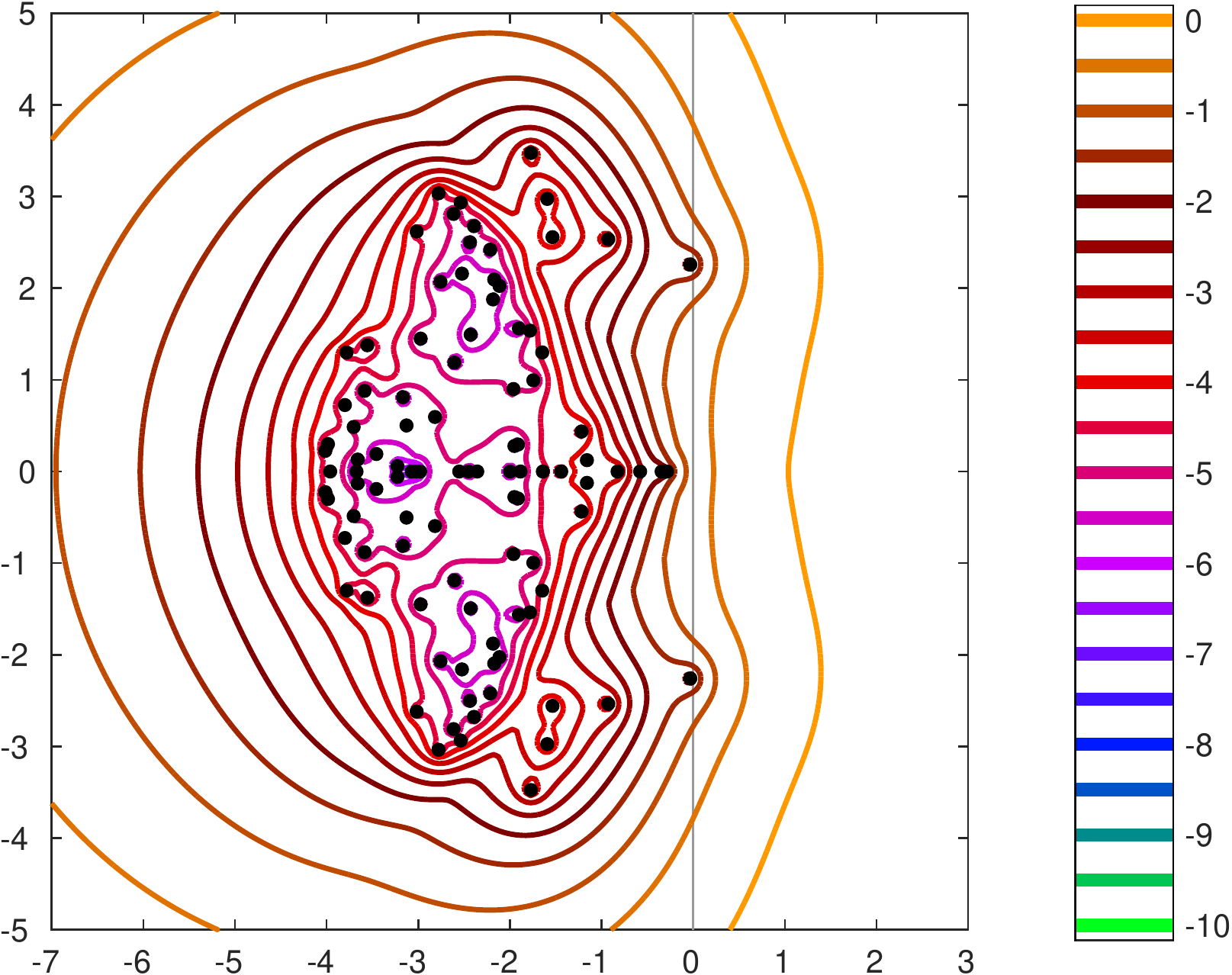}} 
\put(170,175){\includegraphics[scale=0.42]{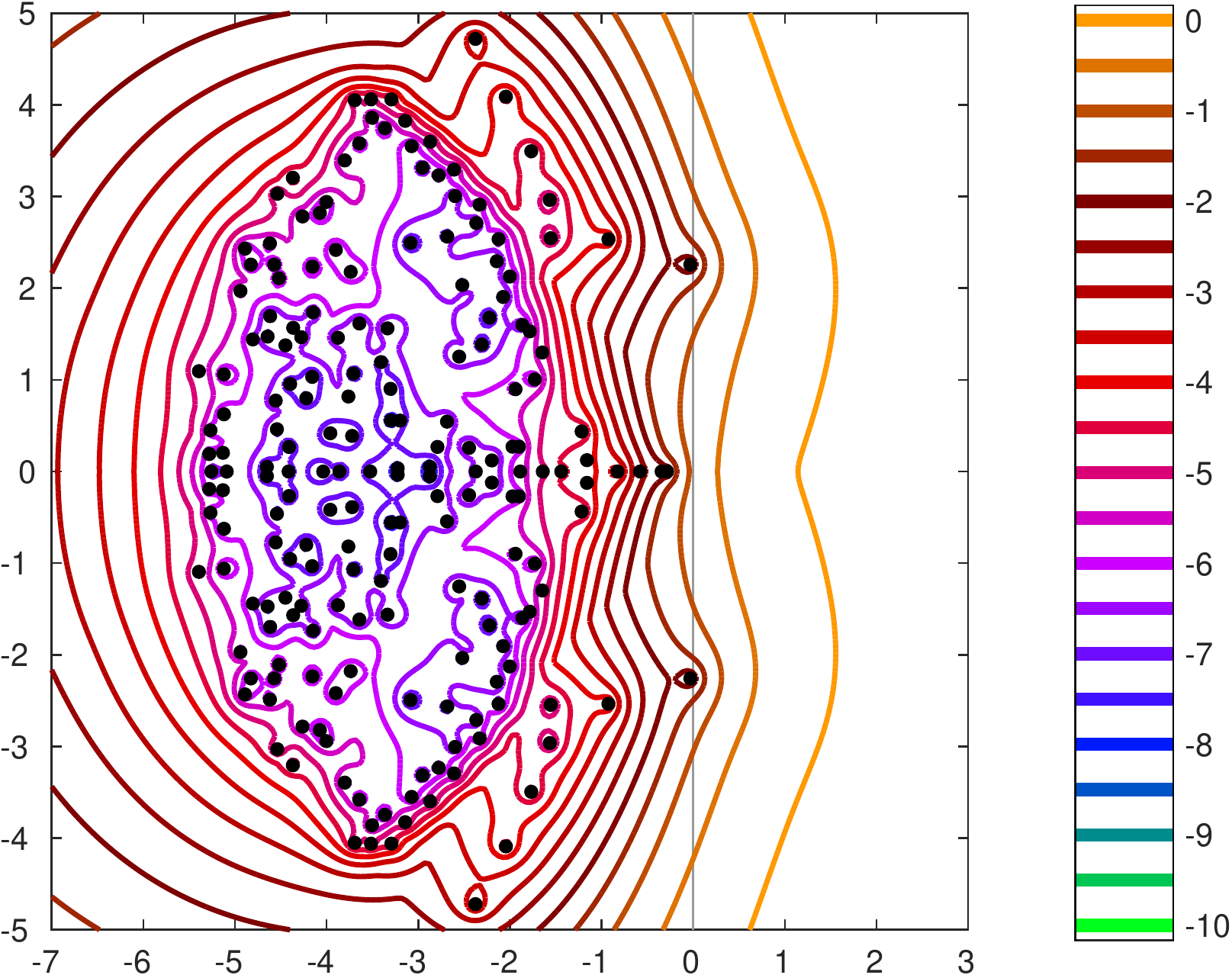}}
\put(0,0){\includegraphics[scale=0.42]{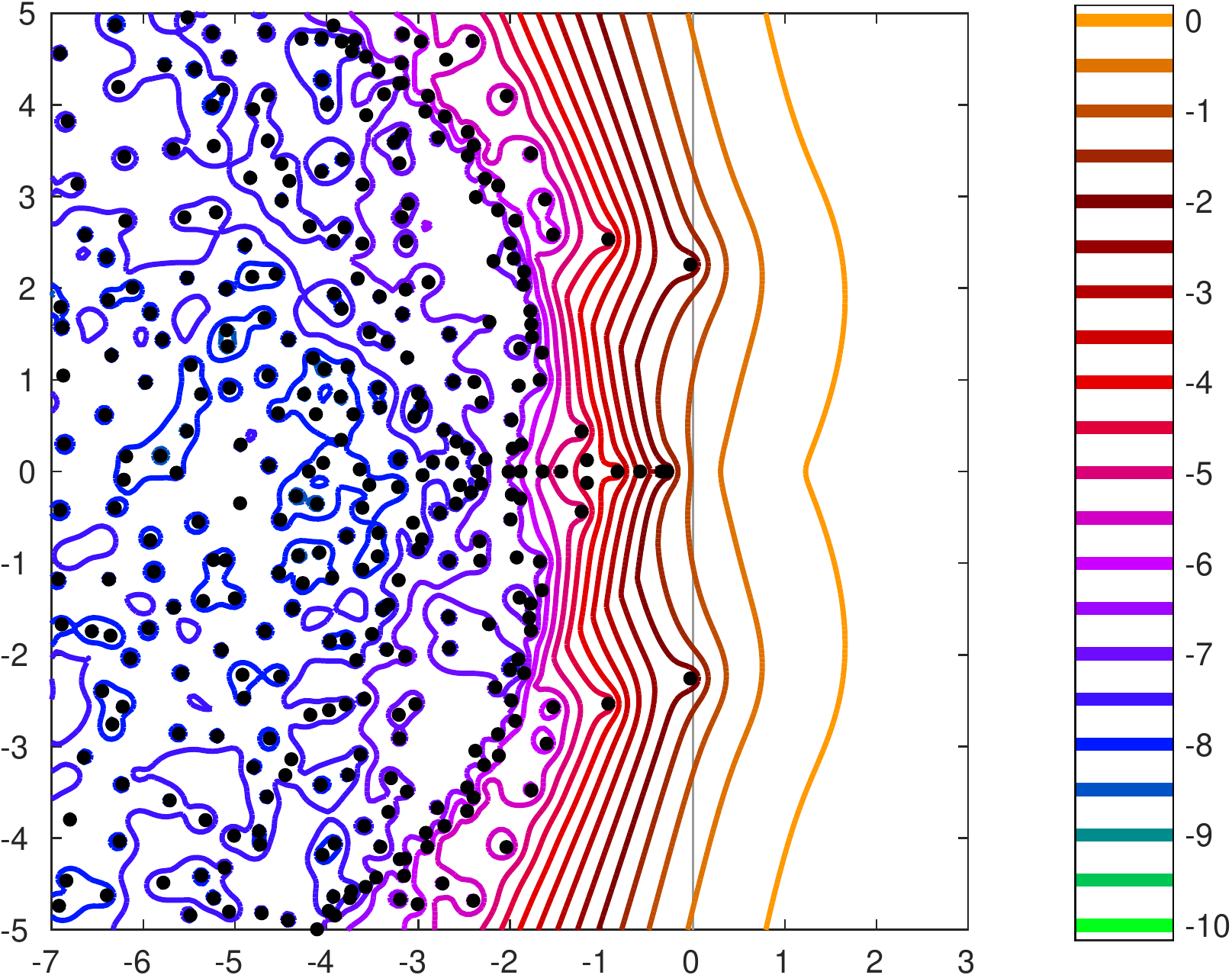}} 
\put(170,0){\includegraphics[scale=0.42]{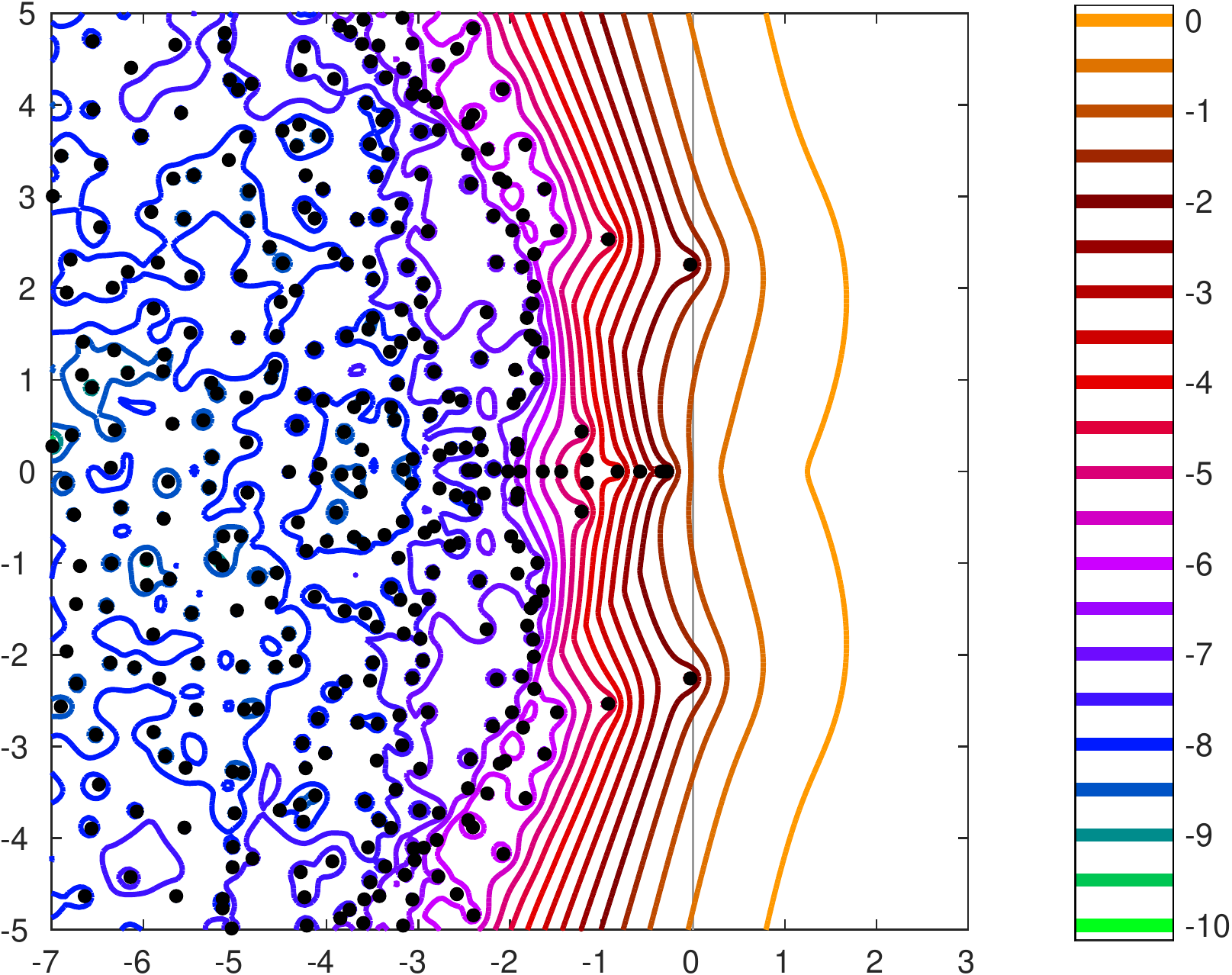}}
\put(10,507){\footnotesize ${\tt nc}=4, k=800$}
\put(180,507){\footnotesize ${\tt nc}=6, k=800$}
\put(10,332){\footnotesize ${\tt nc}=7, k=100$}
\put(180,332){\footnotesize ${\tt nc}=7, k=200$}
\put(10,157){\footnotesize ${\tt nc}=7, k=800$}
\put(180,157){\footnotesize ${\tt nc}=7, k=1600$}
\end{picture}

\caption{\label{fig:obspsa}
Approximations of the pseudospectra $\PSA(\BA,\BE)$ 
for $\eps=10^{0}, 10^{-0.5}, \ldots, 10^{-10}$
for flow around an obstacle with viscosity $\nu=1/175$.
The top two plots use projection onto an invariant subspace of dimension $k=800$
for discretizations of dimension $n=2{,}488$ (${\tt nc}=4$) and 
       $n=37{,}168$ (${\tt nc}=6$).
The bottom four plots project a discretization of size
$n=146{,}912$ (${\tt nc}=7$) onto subspaces of dimension 
$k=100$, $200$, $800$, and $1600$.
The labels on the color bar show $\log_{10}\eps$, 
so, e.g., the orange contour on the right corresponds to
$\eps = 10^0$.
}
\end{figure}

Another wrinkle emerges in these practical computations.
The {\tt eigs} command in MATLAB (which calls the ARPACK software~\cite{LSY98})
returns a basis of eigenvectors that is highly ill-conditioned.  (This is no
surprise, given the significance of the $\eps=10^{-10}$ pseudospectrum in the 
bottom plots of \cref{fig:bfspsa}.)  One could respond to this ill-conditioning
by projecting only onto the dominant component of this subspace, or by generating
a orthonormal basis for all $k$ of the ill-conditioned vectors.  All the
computations shown here use the latter option, projecting
onto an approximate invariant subspace.  
(Some of the interior eigenvalues and 
pseudospectral boundaries for small $\eps$ shown here will be sensitive 
to the basis that {\tt eigs} returns, consistent with the large departure
from normality evident for this example; the same will hold for the
next example shown in \cref{fig:obspsa}.)

\subsection{Flow around an obstacle}

Our second example concerns flow about a square obstacle;
see~\cite[sect.~5.2]{EMSW12} for further details about this example.
(Again we take $\mu = 0.25$.)
As the viscosity decreases, a pair of complex conjugate eigenvalues
crosses the imaginary axis \new{into the right-half plane}
at $\nu \approx 0.00537$~\cite[sect.~5.2]{EMSW12}.
\Cref{fig:obspsa} shows approximations to $\PSA(\BA,\BE)$ for this 
example with viscosity $\nu = 1/175$, just on the stable side of the transition to instability.
(On grid ${\tt nc} = 7$, the spectral abscissa is approximately $-0.0310469$.)
Grid ${\tt nc} = 4$ leaves the problem underresolved, and the rightmost eigenvalue
is real.  
For grids ${\tt nc} = 6$ and $7$, the rightmost eigenvalues form a conjugate pair, 
as expected for this problem~\cite{EMSW12}.
Comparing ${\tt nc}=6$ and ${\tt nc}=7$,
the exterior eigenvalues on the right of the spectrum appear well converged.
For ${\tt nc}=7$, the eigenvalues in the left of the plots change quite a bit as the 
subspace dimension $k$ increases, suggesting that the associated component of the
computed invariant subspace is inaccurate.

\begin{figure}
\begin{center}
%
\includegraphics[scale=0.5]{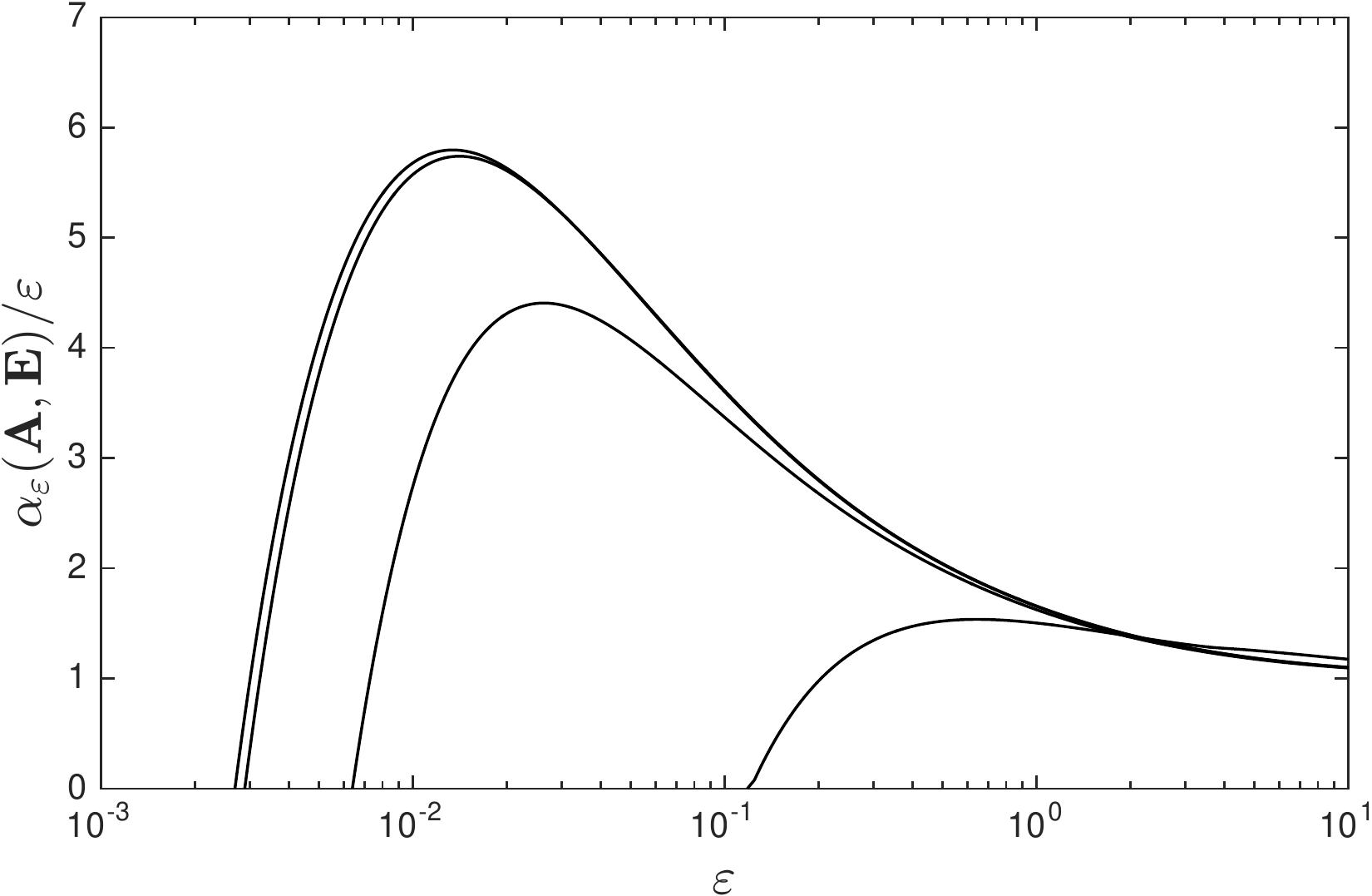}
\begin{picture}(0,0)
\put(-105,135){\textsl{obstacle example}}
\put(-105,122){$k=800$ for all cases}
\put(-201,65){\rotatebox{75}{\footnotesize ${\tt nc} = 7$}}
\put(-192,60){\rotatebox{73}{\footnotesize ${\tt nc} = 6$}}
\put(-174,50){\rotatebox{71}{\footnotesize ${\tt nc} = 5$}}
\put(-105,24){\rotatebox{39}{\footnotesize ${\tt nc} = 4$}}
\end{picture}
\end{center}

\vspace*{-5pt}
\caption{\label{fig:kreiss} Approximations of $\alpha_\eps(\BA,\BE)/\eps$ as a
function of $\eps$ for the obstacle \new{example},
indicating the presence of transient growth. 
The viscosity and projection subspace dimension are fixed
($\nu=1/175$ and projection subspace dimension $k=800$)
while the discretization parameter is varied $({\tt nc} = 4, 5, 6, 7$).}
\end{figure}
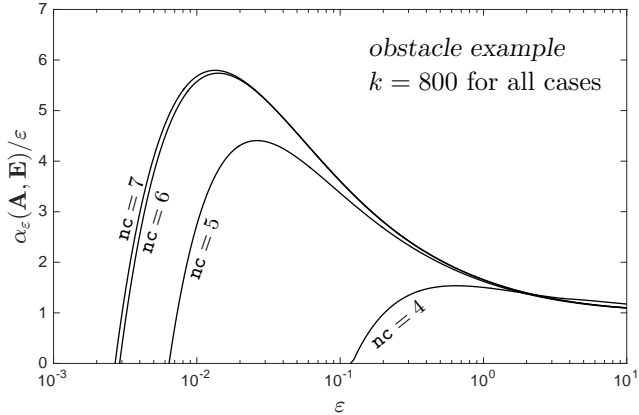

\subsection{Pseudospectral abscissa computations} \label{sec:kreisscomp}
While \cref{fig:bfspsa,fig:obspsa} confirm 
that both flow examples experience transient growth, 
the extent of this growth is difficult to \new{accurately read off}
from plots of the pseudospectra.
\Cref{fig:kreiss} quantifies this growth
by plotting the critical ratio $\alpha_\eps(\BA,\BE)/\eps$ for a range of $\eps$ values
for the obstacle flow problem.
By \Cref{thm:explow}, this ratio provides a \emph{lower bound} on the factor
by which solutions to $\BE\Bx'(t) = \BA\Bx(t)$ can grow.
To make these plots, we used projection onto $k=800$ dimensional invariant
subspaces to estimate $\alpha_\eps(\BA,\BE)$ via~\cref{eq:psnabs}
at hundreds of $\eps$ values 
using the criss-cross algorithm of Burke, Lewis, and Overton~\cite{BLO03b},
as implemented by Mengi, Mitchell, and Overton in EigTool~\cite{Wri02a}.
(One could instead attempt to tackle the large-scale problem directly,
without projection, using alternative algorithms designed to compute the
pseudospectral abscissa of large matrices~\cite{GL11,GO11,KV14}.)

\Cref{fig:kreiss} shows how
$\alpha_\eps(\BA,\BE)$ depends on the quality of the discretization
(for the fixed subspace dimension $k=800$).
For ${\tt nc} = 4$, the plot suggests only mild transient growth; 
larger values of ${\tt nc}$ show more pronounced growth, and appear to be
converging toward a limit: some initial conditions can grow by a factor
of nearly six (at least) before decaying.


\section{Discrete time systems} \label{sec:discrete}
A referee helpfully observed that the analysis described 
above can be adapted to discrete-time systems of the form
\begin{equation} \label{eq:difference}
   \BE@\Bx_{k+1} = \BA@\Bx_k, \rlap{$\qquad \mbox{for $k=0,1,2,\ldots$}$} 
\end{equation}
with initial condition $\Bx_0\in\Cn$.  Using the notation of \cref{sec:definition},
one can write the solution in the form $\Bx_k = \BQ_\mu \By_k + \wt{\BQ}_\mu\Bz_k$
for each $k$.  The difference--algebraic equation~\cref{eq:difference} is equivalent
to 
\begin{align}
    \BG_\mu \By_{k+1} + \BD_\mu \Bz_{k+1} &= (\BI+\mu\BG_\mu)\By_k + \mu\BD_\mu \Bz_k  
     \label{eq:difference1} \\[.25em]
                        \BN_\mu \Bz_{k+1} &= (\BI+\mu\BN_\mu)\Bz_k.
     \label{eq:difference2}
\end{align}
Premultiplying this last equation by $\BN_\mu^{d-1}$ shows that  
$\BN_\mu^{d-1}\Bz_j=\Bzero$ for all $j$.
Now premultiplying~\cref{eq:difference2} by $\BN_\mu^{d-2}$ gives
$\BN_\mu^{d-2}@\Bz_j=\Bzero$ for all $j$.
Repeating this procedure leads, in perfect parallel to the continuous-time case,
to the conclusion that $\Bz_j=\Bzero$ for all $j$;
in particular, $\Bz_0 = \Bzero$ and the initial state must satisfy $\Bx_0\in\Ran(\BQ_\mu)$. 
The system~\cref{eq:difference1} reduces to 
$\BG_\mu\By_{k+1} = (\BI+\mu\BG_\mu)\By_k$, and, 
provided $\Bx_0 \in \Ran(\BQ_\mu)$, equation~\cref{eq:difference}
has the unique solution \[ \Bx_k = \BQ_\mu^{} (\BG_\mu^{-1}+\mu@\BI)^k@\BQ_\mu^*@\Bx_0.\]
The discrete-time iterates thus satisfy
$\|\Bx_k\| \le \|(\BG_\mu^{-1}+\mu@\BI)^k\|@\|\Bx_0\|$, and for each $k$ there
exists some $\Bx_0 \in \Ran(\BQ_\mu)$ for which equality is attained.
Hence, the same $\PSA(\BA,\BE)$ proposed in \cref{def:psa} can be used to bound
the transient behavior, and analogues of many of the theorems in \cref{sec:transient}
follow, with the pseudospectral radius \new{playing a role like} the pseudospectral abscissa \new{in the continuous-time bounds}; 
see~\cite[chap.~16]{TE05} for details.

\section{Conclusions}

What role should structure play in perturbation theory?
This question can be quite delicate, with its answer depending 
on the particular insight one seeks about a given system.
Here we have proposed a definition of the pseudospectrum of 
a matrix pencil that accounts for the structure induced by
a related differential--algebraic equation,
a definition that, by design, gives insight into the transient
dynamics of solutions to the DAE.\ \   
The proposed pseudospectra can be approximated using the 
standard tools for computing rightmost eigenvalues in linear
stability analysis, though the fluid examples shown 
in the last section illustrate that many rightmost eigenvalues
might be required to fully capture the nonnormal dynamics of a
complicated large-scale system.
Further work is needed to understand how inaccuracies in the
computed invariant subspaces affect the approximate pseudospectra, 
and the extent to which reduced order models for descriptor systems 
preserve these pseudospectra and the associated transient dynamics.
The definition here suggests other avenues for investigation, such
as the application of these ideas to DAE systems with polynomial 
structure (as could arise, e.g., from damped mechanical systems with
algebraic constraints), and whether this definition of $\PSA(\BA,\BE)$,
which was motivated by transient analysis rather than eigenvalue perturbations,
might give some insight into the distance of the DAE from instability.

\section*{Acknowledgements}
We thank two referees for their thorough and thoughtful suggestions,
Howard Elman, Rich Lehoucq, Volker Mehr\-mann, and
Paul Van Dooren for helpful discussions about this work, 
and Jonathan Baker for insightful comments on an earlier version on the manuscript.
Tim Davis provided timely advice about the sparse direct solver
in MATLAB, which is heavily used in the numerical examples in \cref{sec:examples}.
We are grateful for support from the Einstein Stiftung Berlin, which enabled the
first author to visit the Technical University of Berlin at a critical point
in this work.

\end{document}